\documentclass[preprint, 3p]{elsarticle}

\pdfoutput=1

\usepackage{lipsum}
\makeatletter
\def\ps@pprintTitle{%
	\let\@oddhead\@empty
	\let\@evenhead\@empty
	\def\@oddfoot{}%
	\let\@evenfoot\@oddfoot}
\makeatother

\usepackage[english]{babel}
\usepackage[utf8]{inputenc}
\usepackage[colorinlistoftodos]{todonotes}
\usepackage{amsmath}
\usepackage{amssymb}
\usepackage{amsthm} 
\usepackage{mathtools}
\usepackage{mathrsfs}
\usepackage{siunitx}
\usepackage{gensymb}
\usepackage{comment}

\usepackage[hidelinks]{hyperref}
\usepackage{algorithm}
\usepackage{algpseudocode}

\usepackage{subfig}
\usepackage{makecell,pict2e}
\usepackage{array, diagbox}
\usepackage{pdflscape}
\usepackage{tikz}

\usepackage{tabularx}
\usepackage{xcolor,colortbl}
\usepackage{multirow}
\usepackage{booktabs}

\usepackage{empheq}

\usepackage[toc,page]{appendix}

\theoremstyle{plain}
\theoremstyle{plain}
\theoremstyle{plain}
\theoremstyle{plain}
\theoremstyle{plain}
\theoremstyle{definition}
\theoremstyle{remark}\newtheorem{remark}{Remark}

\everymath{\displaystyle}

\newif\ifdouble

\newcommand{\papertitle}{An integrated heart-torso electromechanical model for the simulation of electrophysiogical outputs accounting for myocardial deformation}

\newcommand{\keywordOne}{Heart-torso model}
\newcommand{\keywordThree}{Cardiac Electromechanics}
\newcommand{\keywordTwo}{Electrocardiograms}
\newcommand{\keywordFive}{Body Surface Potential Maps}
\newcommand{\keywordFour}{Multiphysics Modeling}

\newcommand{\PhyEP}{\mathscr{E}}

\newcommand{\PhyMec}{\mathscr{M}}

\newcommand{\PhyAct}{\mathscr{A}}

\newcommand{\PhyCirc}{\mathscr{C}}

\newcommand{\PhyCoupl}{\mathscr{V}} 

\newcommand{\VLVthreedim}{V_{\mathrm{LV}}^{\mathrm{3D}}}

\newcommand{\VRVthreedim}{V_{\mathrm{RV}}^{\mathrm{3D}}}

\newcommand{\VLV}{V_{\mathrm{LV}}}

\newcommand{\VRV}{V_{\mathrm{RV}}}

\newcommand{\PLV}{p_{\mathrm{LV}}}

\newcommand{\PRV}{p_{\mathrm{RV}}}

\newcommand{\Circ}{\boldsymbol{c}}

\newcommand{\CircRhs}{\boldsymbol{D}}

\newcommand{\GammaCapsExtLV}{\Gamma_0^{\mathrm{caps,LV}}}
\newcommand{\GammaCapsExtRV}{\Gamma_0^{\mathrm{caps,RV}}}

\newcommand{\GammaCapsLV}{\Sigma_0^{\mathrm{caps,LV}}}
\newcommand{\GammaCapsRV}{\Sigma_0^{\mathrm{caps,RV}}}

\newcommand{\fZero}{{\mathbf{f}_0}}
\newcommand{\sZero}{{\mathbf{s}_0}}
\newcommand{\nZero}{{\mathbf{n}_0}}

\newcommand{\EPchim}{\chi_\mathrm{m}}
\newcommand{\EPCm}{C_\mathrm{m}}
\newcommand{\EPIion}{{\mathcal{I}_{\mathrm{ion}}}}
\newcommand{\EPIapp}{{\mathcal{I}_{\mathrm{app}}}}

\newcommand{\EPdiffTens}{\boldsymbol{D}_i}
\newcommand{\EPdiffTensE}{\boldsymbol{D}_e}

\newcommand{\EPrhsGating}{\boldsymbol{H}}

\newcommand{\Inv}[1]{{\mathcal{I}_{#1}}}
\newcommand{\IIVf}{\Inv{4f}}

\newcommand{\IIVn}{\Inv{4n}}

\newcommand{\displ}{\mathbf{d}} 

\newcommand{\BCmecKepiTens}{\mathbf{K}^{\mathrm{epi}}}
\newcommand{\BCmecCepiTens}{\mathbf{C}^{\mathrm{epi}}}

\newcommand{\GammaBase}{\Gamma_0^{\mathrm{base}}}
\newcommand{\GammaEpi}{\Gamma_0^{\mathrm{epi}}}

\newcommand{\GammaEndoi}{\Sigma_0^{\mathrm{endo,i}}}
\newcommand{\GammaEndoRV}{\Sigma_0^{\mathrm{endo,RV}}}
\newcommand{\GammaEndoLV}{\Sigma_0^{\mathrm{endo,LV}}}

\newcommand{\mecF}{\mathbf{F}}
\newcommand{\mecC}{\mathbf{C}}
\newcommand{\tenspiola}{\mathbf{P}}

\newcommand{\identity}{\mathbf{I}}

\newcommand{\mecNref}{{\mathbf{N}}}

\newcommand{\BCmecVbaseLV}{\mathbf{v}_{\mathrm{LV}}^{\mathrm{base}}}
\newcommand{\BCmecVbaseRV}{\mathbf{v}_{\mathrm{RV}}^{\mathrm{base}}}

\newcommand{\Cai}{{[\mathrm{Ca}^{2+}]_{\mathrm{i}}}}

\newcommand{\SL}{{SL}} 

\newcommand{\Tens}{T_\mathrm{a}} 
\newcommand{\ActStateHF}{\mathbf{s}}

\newcommand{\ActRhs}{\boldsymbol{K}}

\graphicspath{{pictures/}}

\begin{document}
	\begin{frontmatter}
	\title{{\papertitle}}
   
\author[1,2]{Elena Zappon\corref{cor1}}
\ead{elena.zappon@medunigraz.at}
\author[1,3]{Matteo Salvador}
\ead{msalvad@stanford.edu}
\author[1]{Roberto Piersanti}
\ead{roberto.piersanti@polimi.it}
\author[1]{Francesco Regazzoni}
\ead{francesco.regazzoni@polimi.it}
\author[1]{Luca Dede'}
\ead{luca.dede@polimi.it}
\author[4,5]{Alfio Quarteroni}
\ead{alfio.quarteroni@polimi.it}

\cortext[cor1]{Corresponding author}

\address[1]{MOX - Dipartimento di Matematica, Politecnico di Milano, Milan, Italy}
\address[2]{Gottfried Schatz Research Center Biophysics, Medical University of Graz, Graz, Austria}
\address[3]{Institute for Computational and Mathematical Engineering, Stanford University, California, USA}
\address[5]{Politecnico di Milano, Milan, Italy (Professor Emeritus)}
\address[4]{\'Ecole Polytechnique Federale de Lausanne, Lausanne, Switzerland (Professor Emeritus)}

	\begin{abstract}
		When generating in-silico clinical electrophysiological outputs, such as electrocardiograms (ECGs) and body surface potential maps (BSPMs), mathematical models have relied on single physics, i.e. of the cardiac electrophysiology (EP), neglecting the role of the heart motion. Since the heart is the most powerful source of electrical activity in the human body, its motion dynamically shifts the position of the principal electrical sources in the torso, influencing electrical potential distribution and potentially altering the EP outputs. In this work, we propose a computational model for the simulation of ECGs and BSPMs by coupling a cardiac electromechanical model with a model that simulates the propagation of the EP signal in the torso, thanks to a flexible numerical approach, that simulates the torso domain deformation induced by the myocardial displacement. Our model accounts for the major mechano-electrical feedbacks, along with unidirectional displacement and potential couplings from the heart to the surrounding body. For the numerical discretization, we employ a versatile intergrid transfer operator that allows for the use of different Finite Element spaces to be used in the cardiac and torso domains. Our numerical results are obtained on a realistic 3D biventricular-torso geometry, and cover both cases of sinus rhythm and ventricular tachycardia (VT), solving both the  electromechanical-torso model in dynamical domains, and the classical electrophysiology-torso model in static domains. By comparing standard 12-lead ECG and BSPMs, we highlight the non-negligible effects of the myocardial contraction on the EP-outputs, especially in pathological conditions, such as the~VT.

	\end{abstract}
	
	\end{frontmatter}

	\section{Introduction}
Computer-based simulations in cardiac electrophysiology (EP) have significantly advanced over the past decade \cite{trayanova2016,niederer2019computational,corral2020digital,peirlinck2021precision,TRAYANOVA202489}. Among the clinically reproducible outputs of interest, the electrocardiogram (ECG) is a non-invasive and easy-to-achieve recording of cardiac EP, and serves in clinical practice as a default tool. In computational cardiology, the ability of accurately reproduce ECG waveform represents a footprint of quality standard \cite{Gillette2023,zettinig2013,ZETTINIG20141361}. This alignment underscores the potential for simulated models to closely mirror real-world physiological processes \cite{Gillette2023,grandits2023digital,camps2023digital,Qiao2023,GILLETTE2021102080,Okada2013}, including cardiac diseases such as myocardial infarction \cite{Wang2021} and ventricular tachycardia (VT) \cite{lopez2019personalized,Relan2011}. On the other hand, body surface potential maps (BSPMs) are not often utilized in clinical settings, but provide a detail description of bioelectric signals spanning the entire thorax. This allows for a more comprehensive analysis of the electrophysiological condition in contrast to the standard ECG. BSPMs have found application in patients affected by a broad range of pathologies, including myocardial infarction, ventricular hypertrophy, and cardiac arrhythmias \cite{hearts2040040,ISSA2019155}. From a computational perspective, BSPMs have been used in studying ablation targets for atrial arrhythmias \cite{Feng2022,MARQUES2020103904,Ferrer2013} and examining the effects of drug-induced conditions \cite{Zemzemi2013}.

Cardiac EP corresponds to the synchronized depolarization and repolarization of cardiomyocytes.
The tissue depolarization then initiates the mechanical contraction, and therefore deformation, of the heart. Since the heart is the main active source of the electrical field in the human body \cite{Pullan2005,Malmivuo1995}, its deformation leads the potential to (i) influence the transmission of the EP signal within the cardiac tissue itself and (ii) shifts the origins of the electrical field with respect to the body surface, that instead is almost fixed during a single heart beat. Specifically, this last occurrence can lead to changes in the direction of the EP signal propagation throughout the torso, consequently affecting the ECGs and the BSPMs. Indeed, numerous studies have demonstrated that differences in heart location \cite{MACLEOD1998114,Andlauer2018,SCHIJVENAARS2008190,Huiskamp2006HeartPA,NGUYEN2015617,Minchole2019} and dimensions \cite{Andlauer2018,Feldman1985,NAGEL2021102210} in the torso can significantly impact the propagation of the electrical signal within the human body and the EP clinical outputs. This means that the cardiac contraction, determining large deformation and, even, a twist of the heart tissue, may contribute to the genesis of the ECG signals and BSPMs.

A mathematical model for the generation of EP clinical outputs typically involves two key components: a cardiac electrophysiological model, such as the Bidomain \cite{Pullan2005,zappon2023staggeredintime,Poltri2023,franzone2014mathematical,quarteroni2019cardiovascular}, the Monodomain \cite{Pullan2005,franzone2014mathematical,Sundnes2006}, or the Eikonal one \cite{GanderKrauseWeiseretal2023,Stella2022,KONUKOGLU2011134,neic2017efficient},
and a model to calculate the ECG and BSPMs from the electrical signal propagation through the torso \cite{Pullan2005,franzone2014mathematical}.  
The latter can be simulated using the lead-field model \cite{GILLETTE2021102080,multerer2021uncertainty,Potse2018}, enabling the direct computation of ECG leads through time, or by solving a Laplace problem \cite{zappon2023staggeredintime,franzone2014mathematical,Boulakia2010,Potse2009,Aoki1987}, allowing the computation of ECGs and BSPMs as post-processing of the problem solution. These models are usually approximated by means of the Finite Element methods \cite{Pullan2005}. An alternative is represented by the boundary element method (BEM) \cite{sedova2023localization,potse2009cardiac,FUCHS20011400,geselowitz1967bioelectric}.  
Regardless of the chosen representation, the electrophysiology-torso models consider the heart as a static domain - default anatomy is that of the heart in diastasis configuration - embedded in a corresponding static torso domain, therefore neglecting the dynamic cardiac deformations.   
	
The impact of cardiac mechanical displacement on the EP has been explored in some preliminary studies by de Oliveira et al. \cite{deOliveira2013} and Favino et al. \cite{Favino2016}, who employ an electromechanical model (EM) in place of an EP one, yet only simulating electrogram signals (EGMs). In the work of Smith et al. \cite{Smith2003}, EP signals obtained by employing dynamic heart-torso domains depending on cardiac displacement are analyzed, albeit in a 2D context. Conversely, simulations in 3D domains were conducted in studies by Keller et al. \cite{Keller2011}, Wei et al. \cite{Wei2006}, and Xia et al. \cite{Xia2005}, deriving the displacement data from MRI images.
Moreover, different electromechanical modeling approaches yield very different outcomes. Indeed, when accounting for the EGMs in 2D models, Wei et al. \cite{Wei2006} and Xia et al. \cite{Xia2005} observed only a marginal change in torso potential due to cardiac displacement. Smith et al. \cite{Smith2003} noted instead a significant impact of myocardial contraction on the ECGs, particularly evident in the precordial leads and T wave. Keller et al. \cite{Keller2011}, on the other hand, reported a substantial influence of ventricular displacement on lead II. The computation of ECGs and BSPMs in dynamic heart and torso domains has been instead explored only in \cite{moss2021fully}, where the authors incorporate the computed myocardial displacement in the EP outputs simulation by remeshing a portion of the torso for each new cardiac configuration. The procedure of \cite{moss2021fully} is tested only in sinus rhythm conditions. However, in \cite{salvador2022role, salvador2021electromechanical}, the authors emphasize that when pathological conditions and rhythm disorders, such as VT, arise, mechano-electric feedbacks (MEFs) may have a strong effect on the EP signal propagation in the heart and, thus, on the overall EP outputs.

\begin{figure}[!t]
	\centering
	\includegraphics[width=0.8\textwidth]{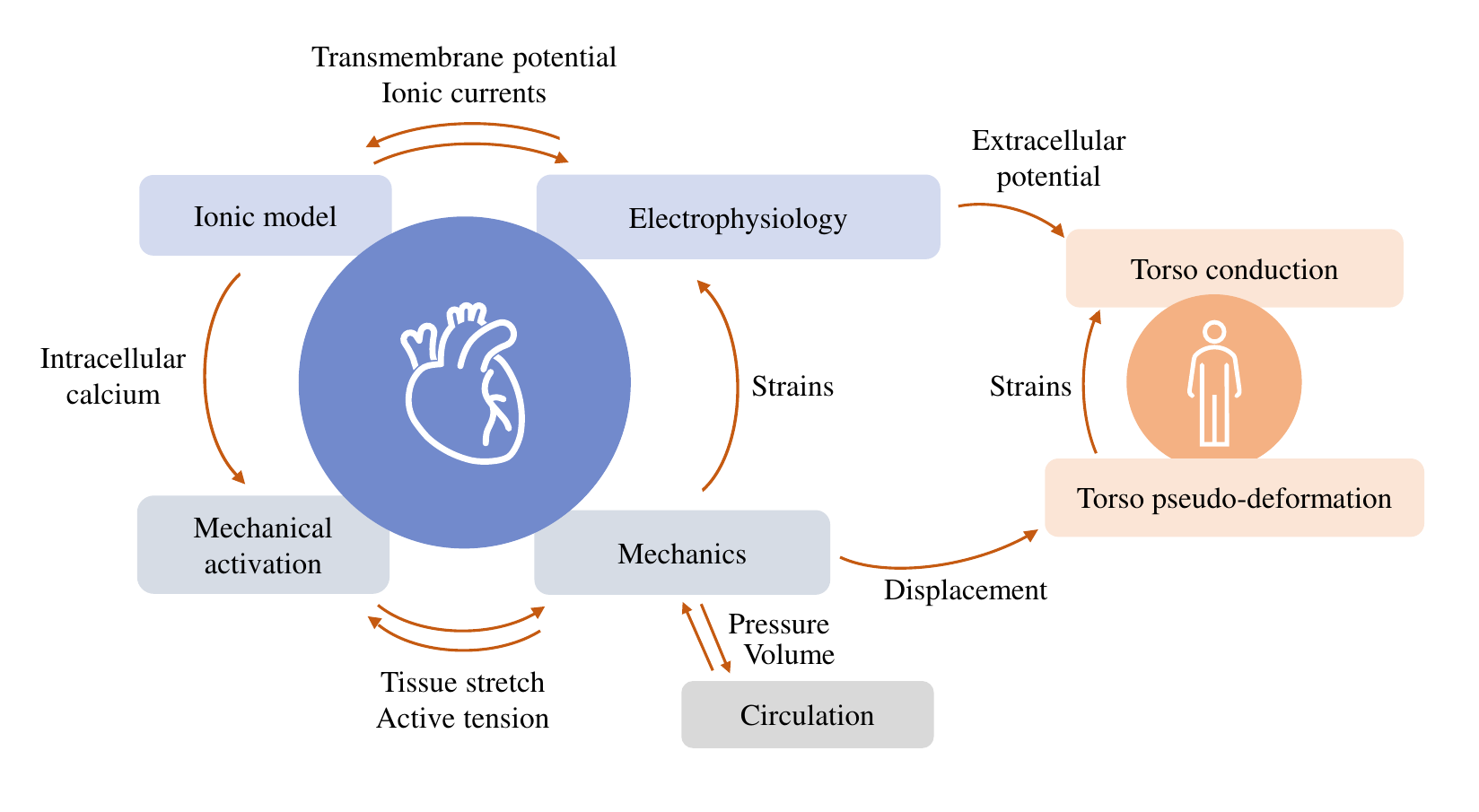}
	\caption{Schematic representation of the implemented electro-mechano-torso model.}
	\label{Fig:EMT_scheme}
\end{figure}

In this work, we propose a 3D multi-physics and multi-scale model for the computation of EP clinical outputs for a moving heart. We couple the cardiac electromechanical model presented in \cite{FEDELE2023115983,piersanti2021closedloop}, with a Laplace problem in the torso \cite{zappon2023staggeredintime}, and employing a lifting technique to account for the dynamic deformation of the torso domain due to the myocardial displacement (see Figure \ref{Fig:EMT_scheme}). The resulting electro-mechanical-torso (EMT) model is an effective alternative of the one presented in \cite{moss2021fully}. The EM model of \cite{FEDELE2023115983,piersanti2021closedloop} includes an accurate description of cardiac EP, extended to represent extracellular potential, passive mechanics, ventricular active contraction, and a reduced-order representation of the circulatory system, fully coupled with the mechanical model of the heart. Moreover we also account for MEFs, fibers-stretch and fibers-stretch-rate feedbacks -- which have been found to be important to regulate the electromechanical behaviour of the heart \cite{FEDELE2023115983}. Finally, the cardiac mechanical displacement is used as boundary conditions to solve a linear elasticity lifting model in the torso domain. The solution of this problem determines the displacement responsible for the online deformation of the torso domain. The Laplace model that accounts for the propagation of the cardiac extracellular potential on the torso \cite{zappon2023staggeredintime} is then solved within the deformed torso domain, and the ECGs and BSPMs are post-processed from the torso solution.  

Concerning the numerical discretization, we extend the accurate and efficient segregated-intergrid-staggered scheme employed in \cite{FEDELE2023115983,REGAZZONI2022111083} to include the torso models. The implemented scheme offers flexibility and allows to prescribe an arbitrary and time-independent displacement to either, or both, cardiac and torso problems.  Consequently, we gain the capability to simulate the EMT model in moving domains while also safely addressing the EMT solution in static domains - thus recovering the static EP-torso model usually employed for computing ECGs and BSPMs. This framework accommodates hybrid scenarios, such as a static heart with a moving torso, or a moving heart with a static torso. This facilitates the separate analysis of (i) the effects of myocardial deformation on EP propagation within the heart and, consequently, on the ECG, and (ii) of the impact of shifting the torso domain shape, and thus the position of the electrical sources in the body, according to the myocardial displacement. The proposed computational framework leverages high-performance computing to enable large-scale simulations, making use of the \texttt{C++} finite element library \texttt{life$^{\text{\texttt{x}}}$} \cite{lifex_2,africa2023lifexep,africa2022lifex}.

We perform numerical simulations using 3D realistic biventricular and torso geometries. To assess the impact of myocardial deformation on the EP outputs, we compare ECGs and BSPMs obtained from the EMT solution in various configurations: dynamic heart-torso domains, static heart-torso domains, and the hybrid configurations. We simulate both healthy and pathological scenarios, specifically replicating VT induced by idealized scar and grey zones on the biventricular septum, following \cite{salvador2022role}.

The work is organized as follows: in Section \ref{Sec:models} we provide the description of EMT mathematical model; in Section \ref{Sec:num_scheme} we briefly describe our numerical framework; in Section \ref{Sec:numerical_results} we present the numerical results obtained with the proposed EMT model; finally, in Sections \ref{Sec:discussion} and \ref{Sec:conclusions} we discuss our key finding and implications of the study, and draw our conclusions.
	\section{Mathematical model}
\label{Sec:models}
\begin{figure}[!t]
	\centering
	\includegraphics[width=0.6\textheight]{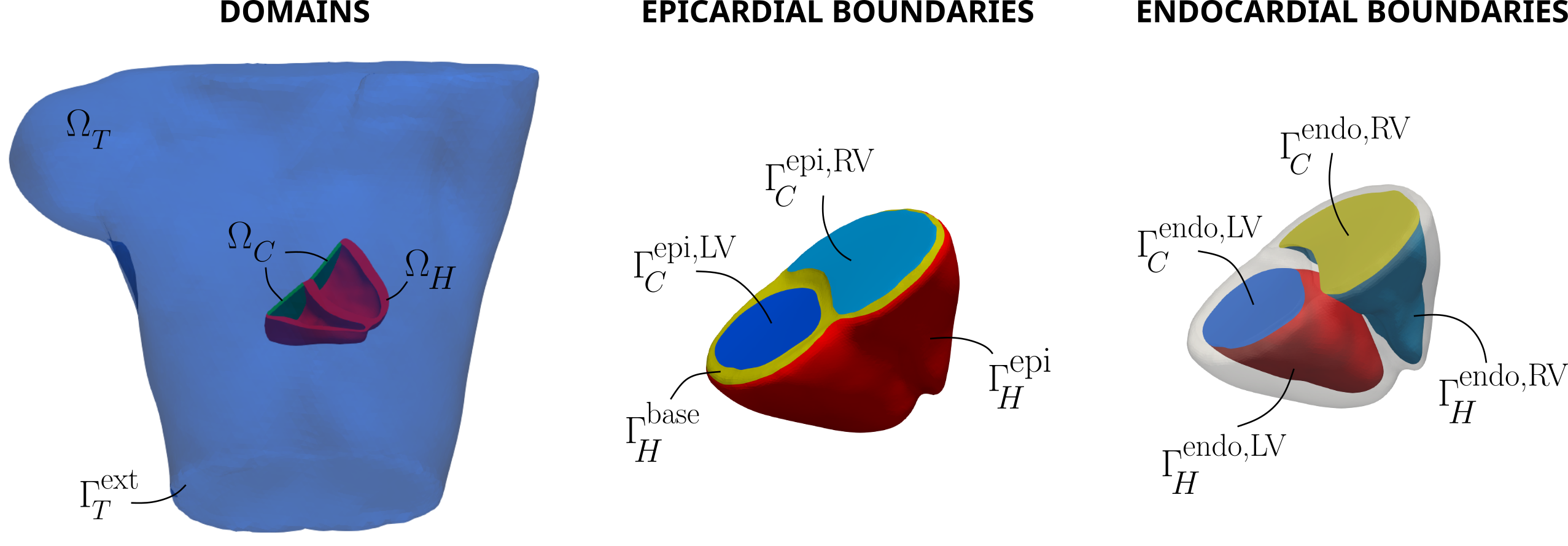}
	\caption{Left: Domains $\Omega_T$ (torso), $\Omega_H$ (biventricular geometry), and $\Omega_C$ (caps). The external surface of the torso is indicated as $\Gamma_T^{\text{ext}}$. Center and right: partitioned of the boundary $\partial \Omega_H$ in epicardium $\Gamma_H^{\text{epi,LV}}$ and $\Gamma_H^{\text{epi,RV}}$, base $\Gamma_H^{\text{base}}$, and left and right endocardium $\Gamma_H^{\text{endo,LV}}$ and $\Gamma_H^{\text{endo,RV}}$ . The portion of $\partial \Omega_C$ representing the surface directed to the torso, both left and right $\Gamma_C^{\text{epi,LV}}$ and $\Gamma_C^{\text{epi,RV}}$, and directed to the cardiac endocardium, both left and right $\Gamma_C^{\text{endo,LV}}$ and $\Gamma_C^{\text{endo,RV}}$, are also indicated.}
	\label{Fig.Domains}
\end{figure}

Let $\Omega_H \subset \mathbb{R}^3$ and $\Omega_T \subset \mathbb{R}^3$ denote two open, time dependent bounded domains, representing the spaces occupied by the human heart and the rest of the human body surrounding $\Omega_H$, respectively, at each time instance of the heartbeat. Hereon, $\Omega_H$ will be represented by a basal biventricular geometry, i.e. a biventricular geometry cut below the cardiac valves (see Figure \ref{Fig.Domains}). To prevent undesirable deformation of the finite elements in $\Omega_T$ around the edges of the cardiac basal plane, and consequent convergence issues of the numerical solver, the internal volume of the ventricular chambers is further isolated from the rest of $\Omega_T$ by sealing the ventricular base with two thin layers of flexible non-conductive tissue, which we refer to as caps and denote by $\Omega_C$ (see Figure \ref{Fig.Domains}). 

\begin{remark}
	While caps do not constitute physiological structures within the heart, in this work they are non-conductive entities from an electrophysiological perspective, and passive tissue from the mechanical standpoint. Furthermore, the applied boundary conditions at their edges ensure synchronized movement with both cardiac tissue and the surrounding torso, without affecting the motion of either heart or torso tissue. Consequently, they do not contribute to the propagation of electrical signals or the deformation of heart and torso domains, thereby having no impact on the generation of ECGs and BSPMs. These caps can be, therefore, entirely disregarded whenever such steep angles are not present in the cardiac domain. 
\end{remark}

The boundary $\partial\Omega_H$ is further split into the left endocardial surface $\Gamma_H^{\text{endo,LV}}$, the right endocardial surface $\Gamma_H^{\text{endo,RV}}$, the epicardial surface $\Gamma_H^{\text{epi}}$, and the base $\Gamma_H^{\text{base}}$. Additionally, we denote:
\begin{itemize} \item $\Gamma_C^{\text{endo,RV}}$ and $\Gamma_C^{\text{endo,LV}}$ as the portions of $\partial \Omega_C$ directed towards the right and left ventricular chambers, respectively.
\item $\Gamma_C^{\text{epi,RV}}$ and $\Gamma_C^{\text{epi,LV}}$ as the portions of $\partial \Omega_C$ directed towards the torso.
\end{itemize}
For the sake of notation, in the rest of the paper we will identify the union of different portions of the boundary of $\Omega_T$, $\Omega_H$ and $\Omega_C$ as:
\begin{itemize}  
	\item $\Gamma = \partial \Omega_T \cap (\partial \Omega_H \cup \partial \Omega_C)$ to represent the interface between the torso domain $\Omega_T$ and the heart-caps volume $\Omega_H$-$\Omega_T$.
	\item $\Gamma^{\text{RV}} = \Gamma_H^{\text{endo,RV}} \cup \Gamma_C^{\text{endo,RV}} \cup \Gamma_C^{\text{epi,RV}}$ for domain boundaries related to the right ventricle.
	\item $\Gamma^{\text{LV}} = \Gamma_H^{\text{endo,LV}} \cup \Gamma_C^{\text{endo,LV}} \cup \Gamma_C^{\text{epi,LV}}$ for domain boundaries related to the left ventricles.
\end{itemize}
Finally, the external surface of the torso is defined as $\Gamma_T^{\text{ext}}$ (see Figure \ref{Fig.Domains}).

For each domain, the deformation is moreover computed starting from a \emph{static reference configuration} denoted by $\Omega^0_{\{H,T,C\}}$.

We denote by $t$ the time variable. To keep the notation light, in this work $t$ is usually omitted. 

The EMT model is obtained by coupling the EM model presented in \cite{FEDELE2023115983,piersanti2021closedloop} (including electrophysiology, active force generation, passive mechanics, and cardiovascular hemodynamics) with  a torso domain deformation model, and a torso passive conduction model. Specifically, the EMT model features the following unknowns:
\begin{equation}
\label{Eq:unknowns}
\begin{array}{lll}
	u : \{\Omega^0_H \cup \Omega^0_C\} \times [0,T] \to \mathbb{R}, && u_e : \{\Omega^0_H \cup \Omega^0_C\} \times [0,T] \to \mathbb{R}, \\
	\boldsymbol{\omega}: \{\Omega^0_H \cup \Omega^0_C\} \times [0,T]  \to \mathbb{R}^{n_{\boldsymbol{\omega}}}, && \mathbf{z}: \Omega^0_H \times[0,T] \to \mathbb{R}^{n_{\mathbf{z}}},\\
	\mathbf{d}_H : \{\Omega^0_H \cup \Omega^0_C\} \times [0,T] \to \mathbb{R}^3, && \mathbf{c} : [0,T] \to \mathbb{R}^{n_{\mathbf{c}}},\\
	p_i: [0,T] \to \mathbb{R}, ~i \in \{RV, LV\}, && \mathbf{d}_T : \Omega_T^0 \times [0,T] \to \mathbb{R}^3, \\
	u_T : \Omega_T^0 \times [0,T] \to \mathbb{R}, && 
\end{array}
\end{equation}
where $u$ and $u_e$ are the transmembrane and extra-cellular potentials, respectively, $\boldsymbol{\omega}$ represents the vector of the ionic variables, $\mathbf{z}$ is the vector of state variables of the force generation model, $\mathbf{d}_H$ is the cardiac mechanical displacement, here extended to include the caps, $\mathbf{c}$ is the state vector of the circulation model, $p_{RV}$ and $p_{LV}$ are the blood pressures inside the biventricular domain, $\mathbf{d}_T$ is the displacement of the torso domain, and $u_T$ is the electric potential in the torso.

In the following sections, we present a brief description of the models involved in the cardiac EM model is presented, and a detailed one for the torso model. A complete description of the EM model can be found in \ref{App:cardiac_EM}.

\subsection{Cardiac electrophysiology}
\label{subsubsec:EP_model}
We simulate the electrical excitation and propagation in the cardiac tissue by solving the Monodomain model \cite{franzone2014mathematical,Pullan2005,Sundnes2006} coupled with the ten Tusscher-Panfilov (TTP06) ionic model \cite{ten2006alternans} in the reference configuration, computing the transmembrane potential $u$ and the ionic variables $\boldsymbol{\omega}$ and $\mathbf{z}$. 

Vectors $\mathbf{f}_0$, $\mathbf{s}_0$, and $\mathbf{n}_0$, denoting the fiber, sheet and sheet-normal directions, are computed using the Laplace-Dirichlet-Rule-Based-Methods (LDRBMs) outlined in \cite{piersanti2021modeling}. A gradual transmural twist from 60$^{\circ}$ at the endocardium to -60$^{\circ}$ at the epicardium is considered for the fiber orientations. Tissue parameters are tuned to achieve conduction velocities of 0.6, 0.4, and 0.2 m/s in the longitudinal $\mathbf{f}_0$, transverse $\mathbf{s}_0$, and normal $\mathbf{n}_0$ directions of cardiac myocytes, respectively \cite{piersanti2021modeling}.

The myocardial conduction system is surrogated by a transmembrane current applied to five activation points on the endocardial surfaces and a thin, fast endocardial layer, able to mimic physiological activation~\cite{durrer1970total,Myerburg1978}. Activation impulses, along with the velocity of the fast endocardial layer, are calibrated to generate physiological ECG waves, as described in \cite{zappon2023staggeredintime}.

To account for the tissue stretch computed by the mechanical model, the deformation gradient tensor $\mathbf{F}_H=\mathbf{I} + \nabla \mathbf{d}_H$ and the corresponding Jacobian $J_H = \text{det}(\mathbf{F}_H) > 0$ are included in the diffusion term of the Monodomain model, as well as in the definition of the anisotropic diffusion tensor (we refer the reader to \cite{piersanti2021closedloop,REGAZZONI2022111083} and \ref{App:cardiac_EM} for a detailed description of the Monodomain formulation).

The extracellular potential $u_e$, which is assumed to be the only electrical potential extending outside the heart and into the rest of the human body \cite{Pullan2005,colli2018numerical}, is obtained from the transmembrane potential $u$ by solving the following Laplace problem \cite{zappon2023staggeredintime,Boulakia2010}:
\begin{subequations}\label{eqn:Pseudo_EM}
	\begin{empheq}[left={\empheqlbrace\,}]{align} 
	-\nabla \cdot(J_H \mecF_H^{-1}(\mathbf{D}_i+  \mathbf{D}_e)\mecF_H^{-T}\nabla u_e) = \nabla \cdot (J_H \mecF_H^{-1}\mathbf{D}_i \mecF_H^{-T}\nabla u)  \qquad \quad  \qquad \:\:\:\, \text{in } \{\Omega_{H}^0\cup \Omega^0_C\}  \times (0,T],\label{eqn:pseudo_EM1}\\
	(J_H \mecF_H^{-1}(\mathbf{D}_i + \mathbf{D}_e)\mecF_H^{-T}\nabla u_e) \mathbf{n}_H = - (J_H \mecF_H^{-1}\mathbf{D}_i \mecF_H^{-T}\nabla u) \mathbf{n}_H  \qquad \quad \:\:\:\, \text{on } \{\partial\Omega_{H}^0\cup \partial\Omega^0_C \} \times (0,T].\label{eqn:pseudo_EM2}
	\end{empheq}
\end{subequations}
In the above equations, $\boldsymbol{D}_{i}$ and $\boldsymbol{D}_{e}$ are the intra-cellular and extra-cellular diffusion tensors, respectively, obtained for each $t \in (o,T])$ as follows:
$$ \boldsymbol{D}_{i,e} = \sigma_{\ell}^{i,e} \frac{\mecF_H\fZero \otimes \mecF_H\fZero}{\|\mecF_H\fZero\|^2} + \sigma_{t}^{i,e} \frac{\mecF_H\sZero \otimes \mecF_H\sZero}{\|\mecF_H\sZero\|^2} + \sigma_{n}^{i,e} \frac{\mecF_H\nZero \otimes \mecF_H\nZero}{\|\mecF_H\nZero\|^2},$$
where $\sigma_{\ell,t,n}^{i,e}$ are the conduction coefficients of the cardiac tissue in the fiber, sheet, and sheet-normal directions. Hereon, the combination of the Monodomain model with problem \eqref{eqn:Pseudo_EM} is referred to as Pseudo-bidomain model.

\begin{remark}
	\label{Rem.Twave}
	In this work, we do not model the effect of transmural and apico-basal cell heterogeneity. Such heterogeneity, that corresponds to variations in action potential durations among cardiac cells, is ultimately responsible for the generation of a physiological T wave. Consequently, the simulated T wave of Section \ref{Sec:numerical_results} may exhibits non-physiological behavior. Despite this limitation, observed variations in the T wave due by myocardial contraction are independent from this heterogeneity, and can be therefore equally appreciated and investigated through this study.
\end{remark}

\subsection{Cardiomyocytes active contraction} 
Cardiomyocyte active contraction bridges electrophysiology and passive mechanics, as it captures the contraction of sarcomeres resulting from fluctuations in calcium concentration $[\text{Ca}^{2+}]_i$, simulated by the ionic model, and the myocardial displacement $\mathbf{d}_H$, provided by the mechanical model. To model the active force generated within the cardiac muscle, we employ the mean-field version of the models proposed in \cite{regazzoni2020biophysically}, denoted as RDQ model. These models have been shown to efficiently describe the regulatory and contractile proteins and their dynamics, ensuring biophysical accuracy \cite{piersanti2021closedloop, fedele2021polygonal}.

The output of the RDQ models is the active tension generated at the microscale, denoted as $T_a$. Furthermore, an active stress approach is used to manage mechanical force generation, receiving intracellular calcium as input from the ionic model. 

\subsection{Passive and active mechanics}
To describe the tissue stress-strain relationship, the cardiac and caps tissue is assumed to be nearly incompressible, anisotropic, and hyperelastic \cite{guccione1991finite,ogden1997non}. The displacement $\mathbf{d}_H$ is therefore obtained by solving a momentum conservation equations endowed with proper boundary conditions. To account for the active mechanics, the Piola-Kirchhoff stress tensor is defined by an orthotropic active stress approach \cite{goktepe2010electromechanics,smith2004multiscale}, with a strain energy term depending on the Usyk constitutive law \cite{guccione1991passive,FEDELE2023115983}. 

The mechanical forces arising from the interaction between the epicardium and the pericardium \cite{Gerbi2018monolithic, pfaller2019importance, strocchi2020simulating} are enforced through generalized Robin boundary conditions applied to the epicardial surface $\Gamma_H^{\text{epi,RV}} \cup \Gamma_H^{\text{epi,LV}}$. Normal stress boundary conditions are employed on the endocardial surfaces and the boundaries of the cap directed towards the torso, that is on $\Gamma^{\text{RV}}$ and $\Gamma^{\text{LV}}$, to model the pressure exerted by the blood. Energy-consistent boundary conditions, addressing the influence of the neglected portion of the heart, i.e. the atria, on the biventricular domain, are applied to $\Gamma_H^{\text{base}}$ \cite{piersanti2021closedloop}. A detailed description of the cardiac mechanical models is given in \ref{App:cardiac_EM}.


\subsection{Circulatory system and coupling conditions}
The influence of the circulatory system on cardiac mechanics is simulated using the 0D closed-loop model proposed in \cite{piersanti2021closedloop,REGAZZONI2022111083}, where other distinct parts of the circulatory system are represented by a series of resistor-inductor-capacitor circuits, the 0D cardiac chambers are characterized by time-varying elastance elements, and heart valves are simulated using non-ideal diodes.

The coupling of the 0D circulatory model with the 3D EM model is achieved by substituting the time-varying elastance elements representing the LV and RV in the circulatory model with their corresponding 3D EM descriptions with suitable volume-consistency coupling conditions \cite{piersanti2021closedloop}. For a more detailed description of the circulation model and the 3D-0D coupling conditions, we refer to \ref{App:cardiac_EM}.

\subsubsection{VT modeling}
VT is a rhythm disorder that is triggered by different cardiac dysfunctions. We define an idealized ischemia in the biventricular geometry in order to create a potential pathway for VT. In this manner, we investigate how the cardiac contraction may influence the generation of pathological ECGs. 

Ischemia is mathematically modeled by defining myocardial regions with slow conduction properties, denoted by gray zones, and actual scars, respectively. To this end, following \cite{salvador2021electromechanical}, we introduce a parameter $\mu = \mu(\mathbf{x}) \in [0,1]$ in the definition of the anisotropic diffusion tensor:
$$ \boldsymbol{D}_{i,e} = \mu\sigma_{\ell}^{i,e} \frac{\mecF_H\fZero \otimes \mecF_H\fZero}{\|\mecF_H\fZero\|^2} + \mu\sigma_{t}^{i,e} \frac{\mecF_H\sZero \otimes \mecF_H\sZero}{\|\mecF_H\sZero\|^2} + \mu\sigma_{n}^{i,e} \frac{\mecF_H\nZero \otimes \mecF_H\nZero}{\|\mecF_H\nZero\|^2},$$
as well as in the TTP06 ionic model \cite{salvador2022role,Arevalo2016}, able to selectively change the conduction properties of specific portion of the cardiac tissue. The parameter $\mu$ is 1 when representing healthy tissue, 0 for scar regions, whereas linear interpolation of $\mu \in [0.1,1]$ can be used to simulate a continuous of gray zones. We refer to \cite{salvador2021electromechanical,salvador2022role} for the description of EM models for ischemic cardiomyopathy. Moreover, since the cardiomyocytes active contraction model receives the intracellular calcium as input from the ionic model, which is dependent of $\mu$, the differentiation between the healthy, scar, and gray zones is straightforward included in the cardiac mechanics.

\subsection{Torso passive conduction and domain deformation}
\label{Sec:Torso}
The computation of ECGs and BSPMs is obtained by solving a Laplace problem that models the torso as a passive conductor, as outlined in \cite{zappon2023staggeredintime,Boulakia2010}, with some modifications. Indeed, the conventional representation views the torso as a static domain. However, the mechanical contractions of the heart result in alterations to the heart-torso interface, affecting the torso domain around the heart. This geometrical deformation changes the amount of tissue through which $u_e$ spreads, potentially influencing the ECG and BSPMs. 
Our model is therefore designed to account for this effect by (i) calculating a virtual deformation of the torso computational domain due to cardiac contraction and (ii) incorporating this deformation into the signal propagation model. 

\begin{remark}
   During each breath, while the heart experiences substantial deformation during the sinus rhythm, the surrounding pericardium remains relatively static \cite{STROCCHI2020109645}. As a result, cardiac contraction minimally influences tissue deformation outside the heart.
   The purpose of our modeling approach is therefore to dynamically deform the domain $\Omega_T$ according to the myocardial contraction, rather than addressing the physiological deformation of the torso tissue around the heart resulting from cardiac displacement. We refer to this non-physiological deformation of the torso as pseudo-deformation, and the corresponding displacement $\mathbf{d}_T$ as pseudo-displacement. 
\end{remark}

The pseudo-deformation of the torso domain caused by myocardial displacement is accounted for through the following linear elasticity problem, as described in \cite{JT_ale, Stein2003}:
\begin{subequations}\label{Eq:lifting_linear_elasticity}
\begin{empheq}[left={\empheqlbrace\,}]{align}
&-\nabla \cdot \boldsymbol{\sigma}(\mathbf{d}_T) = \boldsymbol{0} \qquad \quad \: \text{in }\Omega_T^0, \\
&\mathbf{d}_T = \mathbf{d}_H \qquad \quad  \qquad \quad \:\:\, \text{on }\Gamma,\\
&\mathbf{d}_T = \boldsymbol{0} \qquad \quad \qquad \quad \:\:\:\:\:\ \text{on }\Gamma_T^{\mathrm{ext}}.\label{Eq:still_torso}
\end{empheq}
\end{subequations}
which computes the torso domain pseudo-displacement $\mathbf{d}_T$ induced by the cardiac deformation $\mathbf{d}_H$. In problem \eqref{Eq:lifting_linear_elasticity}, $\boldsymbol{\sigma}$ represents the Cauchy stress tensor:
$$\boldsymbol{\sigma}(\mathbf{d}_T) = \lambda~ tr(\boldsymbol{\varepsilon}(\mathbf{d}_T))\mathbb{I} + 2~\mu~\boldsymbol{\varepsilon}(\mathbf{d}_T),$$
of the torso tissue, being $tr$ the trace operator, $\lambda$ and $\mu$ are the Lamé constants, $\mathbb{I}$ the identity matrix and $\boldsymbol{\varepsilon}$ the strain tensor:
$$ \boldsymbol{\varepsilon}(\mathbf{d}_T) = \frac{1}{2}(\nabla\mathbf{d}_T + (\nabla\mathbf{d}_T)^T).$$
The Lamé constants are expressed in terms of the Young's modulus $E$ and of the Poisson's modulus $\nu$ according to the following formulation:
$$\lambda = \frac{E}{2(1+\nu)}, \quad \quad \mu = \frac{E \nu}{(1+\nu)(1-2\nu)},$$
enabling a direct description of the strain and deformation properties of the material under consideration. 

\begin{remark}
	The external surface of the torso is constrained to be static by Equation \eqref{Eq:still_torso}. Although the torso undergoes physiological motion due to breathing, the frequency of respiratory motion is much lower than that of cardiac contractions \cite{Buehrer2008}. Consequently, we can reasonably assume that fixing $\Gamma_T^{\text{ext}}$ does not result in a loss of information in the simulated ECG and BSPMs.
\end{remark} 

The pseudo-displacement $\mathbf{d}_T$ is utilized to calculate the corresponding gradient deformation tensor $\mathbf{F}_T$ and Jacobian $J_T$. We assume that $\tilde{u}_T$ is the extracellular potential in the torso configuration at a given time $t$, obtained by solving the classical Laplace problem:
\begin{equation}
\begin{cases}
- \nabla \cdot(\mathbf{D}_T\nabla \tilde{u}_T) = 0 &\text{in }\Omega_T(t),\\
(\mathbf{D}_T\nabla \tilde{u}_T) \cdot \mathbf{\tilde{n}}_T = 0 &\text{on }\Gamma_T^{\mathrm{ext}},\\
\tilde{u}_T = \tilde{u}_e &\text{on }\Gamma(t),
\end{cases}
\end{equation} 
where $\mathbf{D}_T$ is the isotropic diffusion tensor in the torso, and is the cardiac extracellular potential at time $t$.
Following the same procedure presented in \cite{colli2016bioelectrical,colli2018numerical} for the  electrophysiological model \eqref{subsubsec:EP_model}, after the pull-back in the reference configuration $\Omega_T^0$, we compute the electrical potential $u_T$ on the reference configuration $\Omega_T^0$ by solving the following problem:
\begin{equation}
\label{Eq:torso_actual_conf}
\begin{cases}
- \nabla \cdot(J_T \mecF_T^{-1}\mathbf{D}_T\mecF_T^{-1}\nabla u_T) = 0 &\text{in }\Omega_T^0,\\
(J_T \mecF_T^{-1}\mathbf{D}_T\mecF_T^{-1}\nabla u_T) \cdot \mathbf{n}_T = 0 &\text{on }\Gamma_T^{\mathrm{ext}},\\
u_T = u_e &\text{on }\Gamma,
\end{cases}
\end{equation}  
where $\mathbf{n}_T$ is the unit outward normal to $\Omega_T^0$.

Problem \eqref{Eq:torso_actual_conf} yields the electric field generated by the heart throughout the entire torso over time. BSPMs are subsequently derived in post-processing as $u_{T_{|\Gamma_T^{\text{ext}}}}$. The traditional leads of the 12-lead ECG system are computed by aggregating the values of $u_T(\mathbf{x}_{e})$ over time, where $\mathbf{x}_{e}$ denotes the spatial position of the electrodes on the body surface (we refer the reader to \ref{App:cardiac_EM} for a detailed description of the leads computation).
Our strategy, by computing the pullback of the moving torso problem in the reference domain, allows for solving the problem on a reference static mesh, thus preventing the need for continuously remeshing the area surrounding the heart, as done e.g. in \cite{moss2021fully}.

\subsection{Reference configurations and initial conditions}
\label{Sub:reference_config}
Cardiac -- and corresponding torso -- geometries are derived from the analysis of in vivo medical images, that are typically acquired during the diastasis phase. The blood pressure acting on the endocardium introduces stress into the resulting geometries. However, the reference configurations $\Omega^0_{H,T,C}$ represent a stress-free state.

While the cardiac geometry undergoes deformation due to active cardiac mechanical dynamics, both the caps and the torso deformations are passively induced by the cardiac displacement. Furthermore, while the entire domain of the caps deforms under the effect of cardiac deformation, the torso domain is altered only at the heart-torso interface. This same dynamic interaction among the three domains is considered in the computation of their respective reference configurations.

Assuming that the imaging of the cardiac geometry $\tilde{\Omega}_H$ refers to a diastasis configuration, the cardiac reference configuration is obtained by solving the inverse problem presented in \cite{REGAZZONI2022111083}. As the motion of the caps depends on the blood pressure in the cardiac chambers, the inverse problem of \cite{REGAZZONI2022111083} is solved in a unique domain $\tilde{\Omega}_H \cup \tilde{\Omega}_C$, constructed by incorporating the artificial caps $\tilde{\Omega}_C$ into the cardiac imaging configuration~$\tilde{\Omega}_H$.

Once the reference configurations $\Omega_H^0$ and $\Omega_C^0$ have been computed, the torso reference configuration $\Omega_T^0$ is constructed by inserting $\Omega_H^0 \cup \Omega_C^0$ into the imaging torso domain $\tilde{\Omega}_T$, replacing $\tilde{\Omega}_H$.	

Initial conditions for the electromechanical problem are set as in \cite{piersanti2021closedloop}, by surrogating the 3D-0D EM problem with a 0D emulator able to efficiently compute the model parameters at limit cycle \cite{regazzoni2021emulator}.
	\section{Numerical framework}
\label{Sec:num_scheme}
The discretization strategy employed for numerically solving the EMT problem is derived by extending the segregated-intergrid-staggered numerical scheme presented in \cite{piersanti2021closedloop} to incorporate the pseudo-deformation of the torso domain and passive conduction. This scheme allows to combine the static and/or dynamic configurations for both the heart and the torso, thereby distinguishing between the effects of cardiac motion and torso pseudo-deformation on the ECG and BSPMs.

Details regarding time and space discretization strategies are provided in Section \ref{Sub:time_space_discr}, while the possible combinations of cardiac and torso domain configurations are outlined in Section \ref{Sub:comb_conf}.

\begin{figure}[!t]
	\centering
	\includegraphics[width=0.75\textheight]{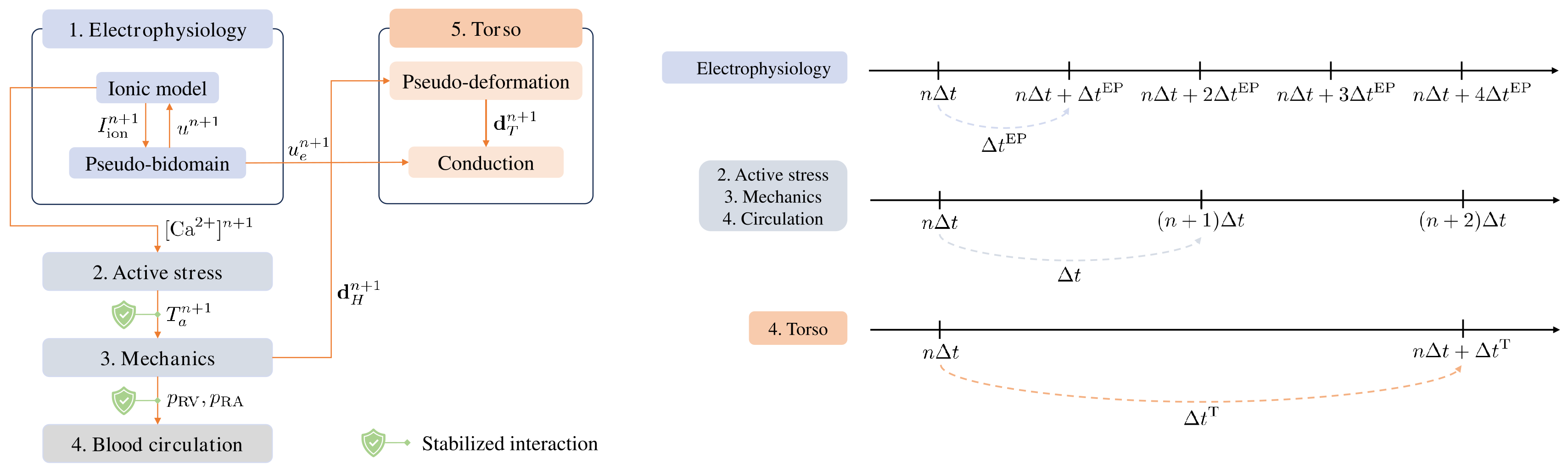}
	\caption{Time-advancing scheme for the coupled EMT model. The number referred to the computational order for a single time step $t_n$.}
	\label{Fig.numerical_scheme}
\end{figure}

\subsection{Time and space discretizations}
\label{Sub:time_space_discr}
The time interval $[0,T]$ is divided into sub-intervals $[t_i, t_{i+1}]$ for $i=0, 1, \dots, N$, where $t_0 = 0$ and $t_N = T$, ensuring a uniform time step size of $\Delta t$ for all $i$. Hereon, we denote the approximation of the solution variables \eqref{Eq:unknowns} at a given time step $t_n$ by the superscript $n$, e.g. $u(t_n) \approx u^n$.

The staggered scheme presented in \cite{piersanti2021closedloop} to solve the coupled EM model is extended to include the pseudo-deformation of the torso domain \eqref{Eq:lifting_linear_elasticity} and the propagation of the cardiac extracellular potential $u_e$ in the torso domain \eqref{Eq:torso_actual_conf} as final step of the simulation, as shown in Figure \ref{Fig.numerical_scheme}.

The scheme involves solving different subproblems separately, using explicit coupling when stability is not a concern. This allows for choosing different time steps for the cardiac EP, mechanics, and torso models. The reference timestep $\Delta t$ is used for the solution of the cardiomyocytes active force generation and the cardiac mechanics. Electrophysiology has a faster dynamics with respect to the other models, requiring a smaller timestep to ensure an accurate solution. Therefore, a timestep $\Delta t^{\text{EP}} = \Delta t/N_{\text{EP}}$, with $N_{\text{EP}} \in \mathbb{N}$, is introduced. Cardiac electrical outputs, such as the ECGs and the BSPMs, can be computed less frequently. Indeed, a third timestep $\Delta t^{\text{T}} = N_{\text{T}} \Delta t$, with $N_{\text{T}} \in \mathbb{N}$, is introduced for this purpose.

Time derivatives are approximated using Finite Difference schemes \cite{quarteroni2010numerical}. The cardiac electrophysiology model is solved using the Backward Differentiation Formula of order 2 (BDF2). An implicit-explicit (IMEX) scheme is adopted to treat the diffusion term implicitly and the ionic and reaction terms explicitly \cite{FEDELE2023115983, niederer2011simulating}. Mechanical activation and mechanics are solved with a BDF1 scheme, employing an explicit method for the cardiomyocyte contraction problem and a fully implicit scheme for the cardiac mechanical model. The circulation model is finally solved explicitly with a Runge-Kutta method of order 4 \cite{FEDELE2023115983}.

\begin{remark}
	The torso problem presented in Section \ref{Sec:Torso} is formally time-independent. Since the boundary conditions of \eqref{Eq:lifting_linear_elasticity} and \eqref{Eq:torso_actual_conf} depend on cardiac quantities varying over time, the solution of the torso problems inherits time-dependency from the cardiac EM models. However, from a numerical standpoint, no time discretization strategies are required for the torso problems.
\end{remark}

The three domains $\Omega_H^0$, $\Omega_C^0$, and $\Omega_T^0$ are spatially discretized on tetrahedral meshes conforming at the domain interfaces. The electrophysiology model is discretized in space using piecewise quadratic finite elements $\mathbb{P}_2$ \cite{Hughes2010, quarteroni2009numerical}, which have been shown in \cite{AFRICA2023111984} to provide improved accuracy with a lower number of degrees of freedom compared to $\mathbb{P}_1$ linear elements for cardiac electrophysiology. 
Linear finite elements $\mathbb{P}_1$ are used for the mechanical activation problem, active and passive mechanics, and torso problems. An efficient intergrid transfer operator is finally used to evaluate the feedback between EP and other physics.

\subsection{Combining static and dynamic heart and torso configurations}
\label{Sub:comb_conf}

Our numerical framework provides not only flexibility in terms of spatial and time discretizations, but allows the exploration of various electrophysiological and torso domain combined configurations. The comprehensive representation of the cardiac EM phenomena encompasses both the dynamic displacement of the heart and the torso domain.

Our numerical scheme facilitates this analysis by allowing the prescription of an arbitrary, time-independent displacements $\mathbf{d}_H$. This static displacement can be utilized:
\begin{itemize}
\item to compute a non-varying deformation tensor $\mathbf{F}_H$ for the Pseudo-bidomain model, simulating cardiac EP on a static configuration, which corresponds to the simulation of the cardiac electrophysiology without mechanical feedback.
\item as a boundary condition for the linear elasticity problem \eqref{Eq:lifting_linear_elasticity} describing torso pseudo-deformation, thus equivalently simulating the torso as a passive conductor in a static domain. 
\end{itemize}
Employing a time-independent $\mathbf{d}_H$ also corresponds to refraining from using the MEFs on either the cardiac or torso domain.
This enables simulation of the EMT model under various scenarios on: 
\begin{enumerate}
	\item both cardiac and torso moving domains, i.e. the most complete representation
	\item both cardiac and torso static domains, recovering the static EP-torso model typically used for computing ECGs and BSPMs
	\item a static heart with a moving torso
	\item a moving heart with a static torso. 
\end{enumerate}
The hybrid configurations 3.~and 4.~allow the analysis of (i) the effects of myocardial deformation on EP propagation within the heart and, consequently, on the ECG and BSPMs (equivalently reading table \ref{Tab.table_configurations_EMT} by columns), and (ii) the impact of shifting the heart-torso interface and thus varying the torso domain according to the myocardial displacement (equivalently reading table \ref{Tab.table_configurations_EMT} by rows).

\begin{table}[!h]
\begin{center}\setlength{\extrarowheight}{2pt}
	\begin{tabular}{c||c|c}
		\hline
		\diagbox[innerwidth = 4cm, height = 6ex]{\textbf{Heart}}{\textbf{Torso}} &Static configuration &Moving configuration \\\hline\hline
		Static configuration & \diagbox[innerwidth = 4cm, height = 8ex, linewidth=1.2pt,linecolor=blue]{No MEFs}{No MEFs} &\diagbox[innerwidth = 4cm, height = 8ex, linewidth=1.2pt,linecolor=red]{No MEFs}{MEFs}\\ \hline
		Moving configuration &\diagbox[innerwidth = 4cm, height = 8ex,linewidth=1.2pt,linecolor=green]{MEFs}{No MEFs} &\diagbox[innerwidth = 4cm, height = 8ex,linewidth=1.2pt]{MEFs}{MEFs}\\ \hline	
 \end{tabular}
\end{center}
\caption{Sketch of the cardiac and torso domain configurations used in the EMT simulations. The two domains can be either in a static configuration -- by prescribing a time-independent $\mathbf{d}_H$ - or in a moving configuration -- dynamically computing $\mathbf{d}_H$. The colors of the diagonal lines refer to the ECG traces of Section \ref{Sec:numerical_results}.}
\label{Tab.table_configurations_EMT}
\end{table}

	\section{Numerical results}
\label{Sec:numerical_results}
We present the results obtained using our EMT model and numerical framework are presented and discussed. The EMT model was use to analyze two scenarios: the first representing healthy conditions, and the second simulating VT. 

The organization of this section is as follows: the volumetric models and baseline parameters common to both test cases are presented in Section \ref{Sub.volum_models}, the baseline simulations representing a healthy patient are displayed in Section \ref{Sub:healthy_cond}, and the results for VT simulations are shown in Section \ref{Sub:VT_conditions}.

\subsection{Numerical simulation settings}
\label{Sub.volum_models}

\begin{figure}[!t]
	\centering
	\includegraphics[width=0.75\textheight]{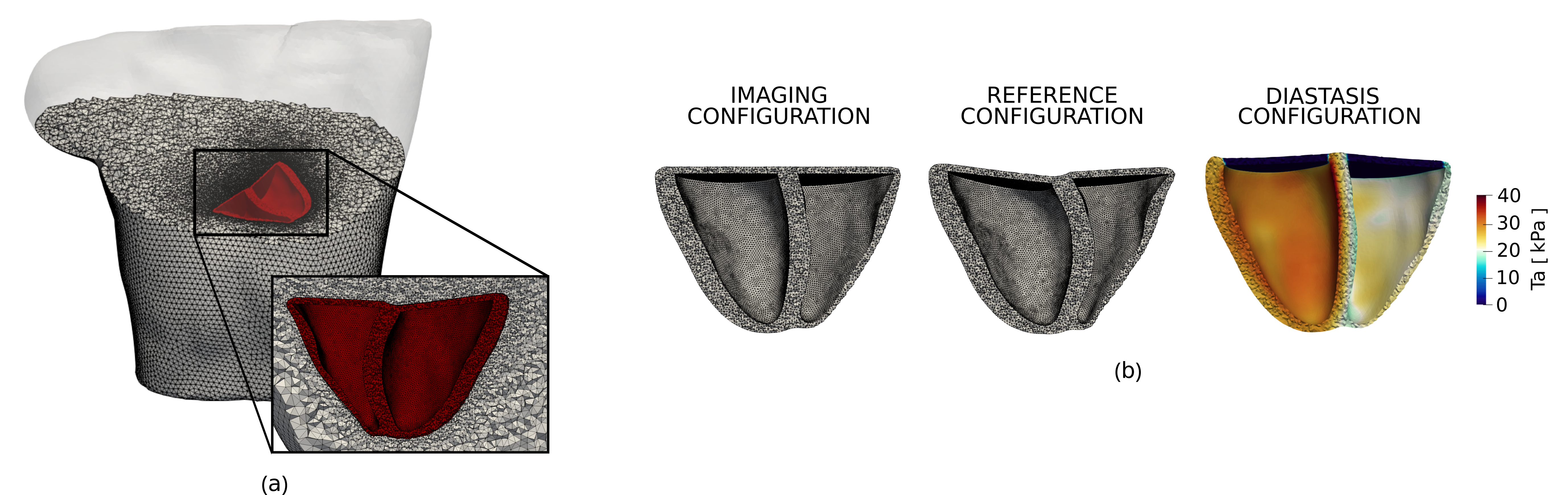}
	\caption{(a) Volumetric geometries for the EMT model. (b) Cardiac imaging configuration, reference configuration $\Omega_H^0$ and diastasis configuration computed from the EM simulation at the diastasis phase of the heart beat. The active tension $T_a$ for the diastasis configuration is also displayed.}
	\label{Fig.mesh_configurations}
\end{figure}

The cardiac biventricular model is based on the Zygote Solid 3D Heart Model \cite{zygote2014}, which is an anatomically accurate CAD model of the entire human heart reconstructed from high-resolution CT scans. It represents a healthy male subject from the 50th percentile of the United States population. The original model has been processed to fit the domain described in Section \ref{Sec:models}. The biventricular geometry is cut below the valves to generate a basal biventricular geometry $\tilde{\Omega}_H$. Two thin layers of tissue are included to close the ventricular chambers at the level of the basal plane, representing $\tilde{\Omega}_C$. The geometrical surfaces are then labeled according to Figure \ref{Fig.Domains}.
Once the tetrahedral mesh is generated, the cardiac reference configuration $\Omega_H^0 \cup \Omega_C^0$ is computed by preprocessing the heart-caps domain $\tilde{\Omega}_H \cup \tilde{\Omega}_C$ as described in Section \ref{Sub:reference_config}. Since the elements in the heart and caps are deformed during the computation of the cardiac reference domains, a final remeshing of the cardiac reference configuration is carried out to improve mesh quality in $\Omega_H^0 \cup \Omega_C^0$ (see Figure \ref{Fig.mesh_configurations}).

The torso volumetric model is derived from the 3D RIUNET torso model \cite{Ferrer2015}, which is publicly available in the online repository of the Center for Integrative Biomedical Computing \cite{CIBC}. The organs, including blood pools, and the atria are removed. The ventricles are substituted with our cardiac reference configuration. The torso reference configuration $\Omega_T^0$ is obtained by generating a tetrahedral mesh within the torso, conforming to $\Omega_H^0 \cup \Omega_C^0$ (referred to Figure \ref{Fig.mesh_configurations} for the torso reference configuration mesh).
All the processing and meshing procedures are performed using the open-source softwares \texttt{vmtk} \cite{fedele2021polygonal,antiga2008image} and \texttt{Paraview} \cite{paraview}.

The cardiac mesh $\Omega_H^0 \cup \Omega_C^0$ is composed of 135K vertices and 715K elements, with an average edge length of \SI{1.6}{\milli\metre}. 
The torso mesh is comprised of an additional 406K vertices and 2.46M elements, with an average edge length of \SI{5.4}{\milli\metre}. 

Both the mechanical and torso simulations are set with the same time step $\Delta t = \Delta t^{\text{T}} =$ \SI{1}{\milli\second}, while  $\Delta t^{\text{EP}} =$ \SI{0.5}{\milli\second} is chosen for the electrophysiological problem, i.e. $N_{\text{EP}} = 20$ and $N_{\text{T}} = 1$. The ionic model is initialized by conducting a 1000-cycle long single-cell simulation for the TTP06 model. Initial circulation variables are calibrated at the limit cycle through the procedure described in Section \ref{Sub:reference_config}. Although five heartbeats are simulated, the results are presented only for on the last two heartbeats to reduce the effects of initialization. 

\begin{figure}[!t]
	\centering
	\includegraphics[width=0.7\textheight]{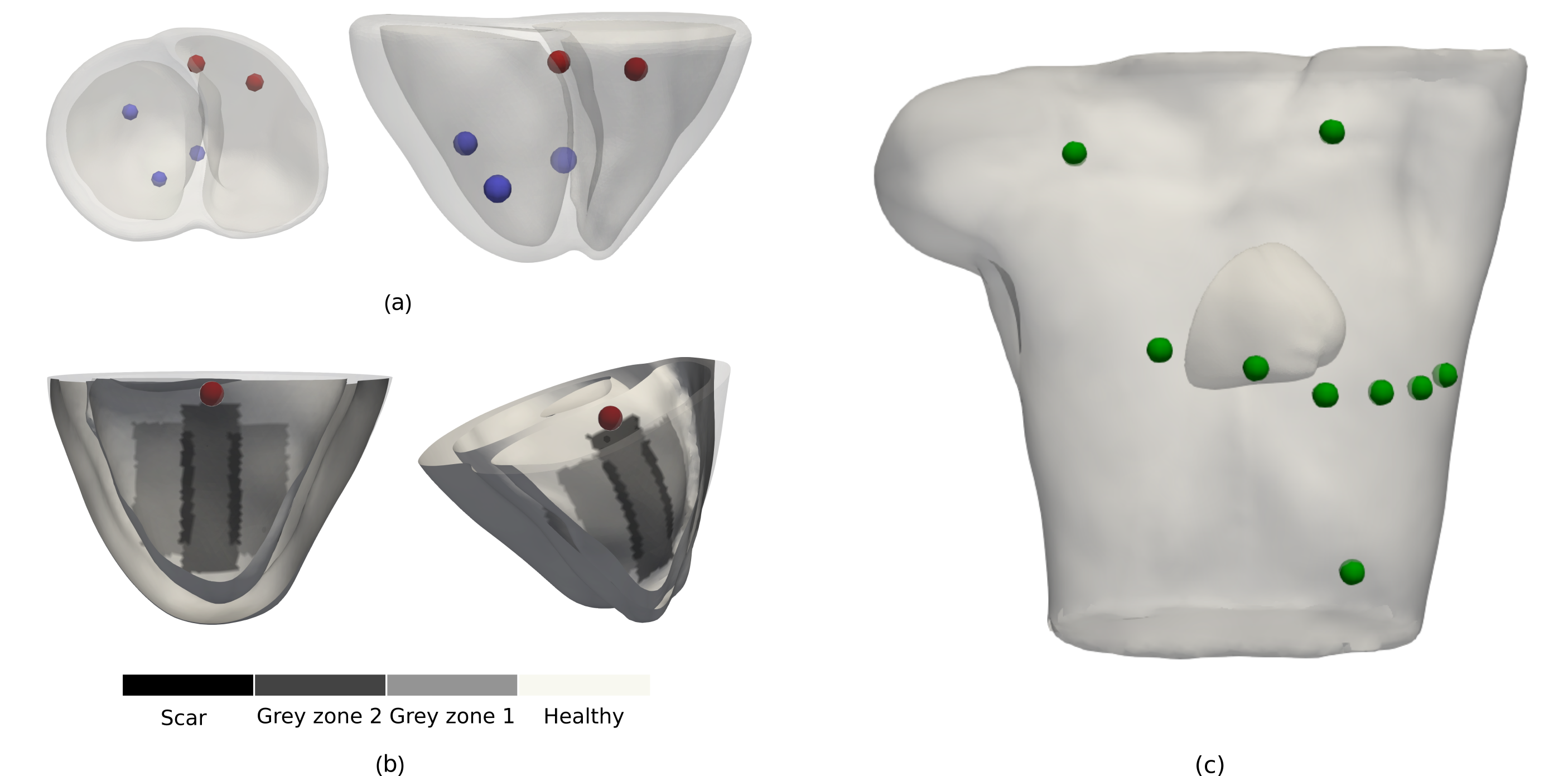}
	\caption{(a) Location of the spherical impulses on the reference cardiac geometry $\Omega_H^0 \cup \Omega_C^0$ used to activate the heart for the test case in healthy conditions. (b) Cardiac reference geometry $\Omega_H^0$ with the idealized distribution of scar (black), grey zones (grey), and healthy tissue (white) over the myocardium. The site of the activation point for the S1, S2, S3 and S4 pacing protocol is also depicted with a red sphere. (c) Position of the 9 electrodes to compute the 12-lead ECG.}
	\label{Fig.activation_numerics}
\end{figure}

In Section \ref{Sub:healthy_cond}, the biventricular conduction system is simplified down to five focal spherical activation points, illustrated in Figure \ref{Fig.activation_numerics}(a), along with a thin, fast conduction layer on the endocardia, following the approach in \cite{zappon2023staggeredintime}. 

In Section \ref{Sub:VT_conditions}, we induced sustained VT with a figure-of-eight pattern through an isthmus, bordered laterally by scars that act as conduction blocks, as shown in Figure \ref{Fig.activation_numerics}(b). Following \cite{salvador2022role, FRONTERA20201719}, we applied an S1-S2-S3-S4 stimulation protocol consisting of four Gaussian stimuli located on the septum near the scar zone (see Figure \ref{Fig.activation_numerics}). The first stimulus S1 is applied at time $t = $\SI{0}{\second}, the second one S2 at $t =$ \SI{0.45}{\second}, the third one S3 at $t =$ \SI{0.75}{\second}, and the fourth one S4 at $t =$ \SI{1.02}{\second}.

The classical 12-leads ECG is computed by placing the nine electrodes in realistic locations on the surface of the torso domain, as depicted in Figure \ref{Fig.activation_numerics}. The same 9 position are used for all the simulations.

In both test cases, all four combinations of cardiac and torso configurations from Section \ref{Sub:comb_conf} are simulated. The static cardiac configuration is represented by the volume occupied by the heart in the diastasis phase (\SI{150}{\milli\second} into the heartbeat), corresponding to classical imaging configurations. However, the original volumetric configuration $\tilde{\Omega}_H$, which ideally represents the heart in the diastasis phase, is never exactly captured in the cardiac electromechanical simulation.
Given the significant role of the cardiac shape in the analysis presented in this work, in order to maximize consistency, the cardiac diastasis configuration is extracted in terms of cardiac displacement $\mathbf{d}_H$ directly from the electromechanical simulation (refer to Figure \ref{Fig.mesh_configurations} for a visual comparison between the reference and diastasis configurations).

The presented numerical framework has been implemented in \texttt{life}$^{\text{x}}$ \cite{lifex_2,AFRICA2022101252}, an in-house high-performance C++ finite element library, based on \texttt{deal.II} \cite{dealII91}, specifically designed for cardiac applications. All numerical simulations were performed using the \texttt{iHeart} cluster, a Lenovo SR950 with 192-Core Intel Xeon Platinum 8160, operating at 2100 GHz, and equipped with 1.7TB RAM, located at MOX, Dipartimento di Matematica, Politecnico di Milano. 
A simulation involving 5 heart-beats is approximately 11 hours long using 48 cores, with only around 6\% dedicated to torso computation. Consequently, the overall cost of an EMT model aligns with that of the EM model \cite{piersanti2021closedloop}.
\subsection{Healthy scenario}
\label{Sub:healthy_cond}

\begin{figure}[!h]
	\centering
	\includegraphics[width=0.75\textheight]{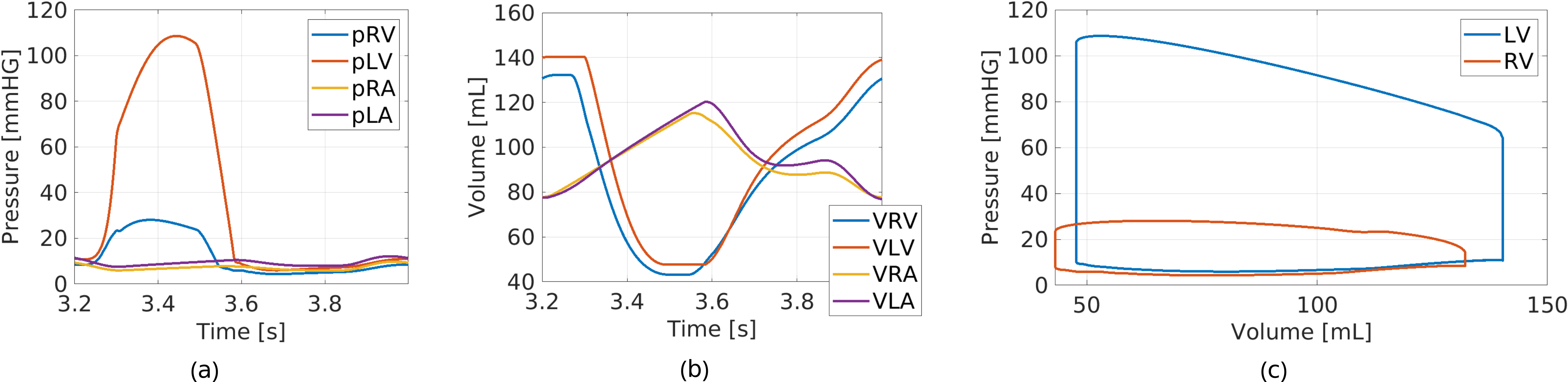}
	\caption{Resulting traces of the circulatory system. (a) Pressures over time. (b) Volumes over time. (c) Pressure-volume loops in the ventricles.}
	\label{Fig.Mechanics_results_loops}
\end{figure}

\begin{figure}[!t]
	\centering
	\includegraphics[width=0.73\textheight]{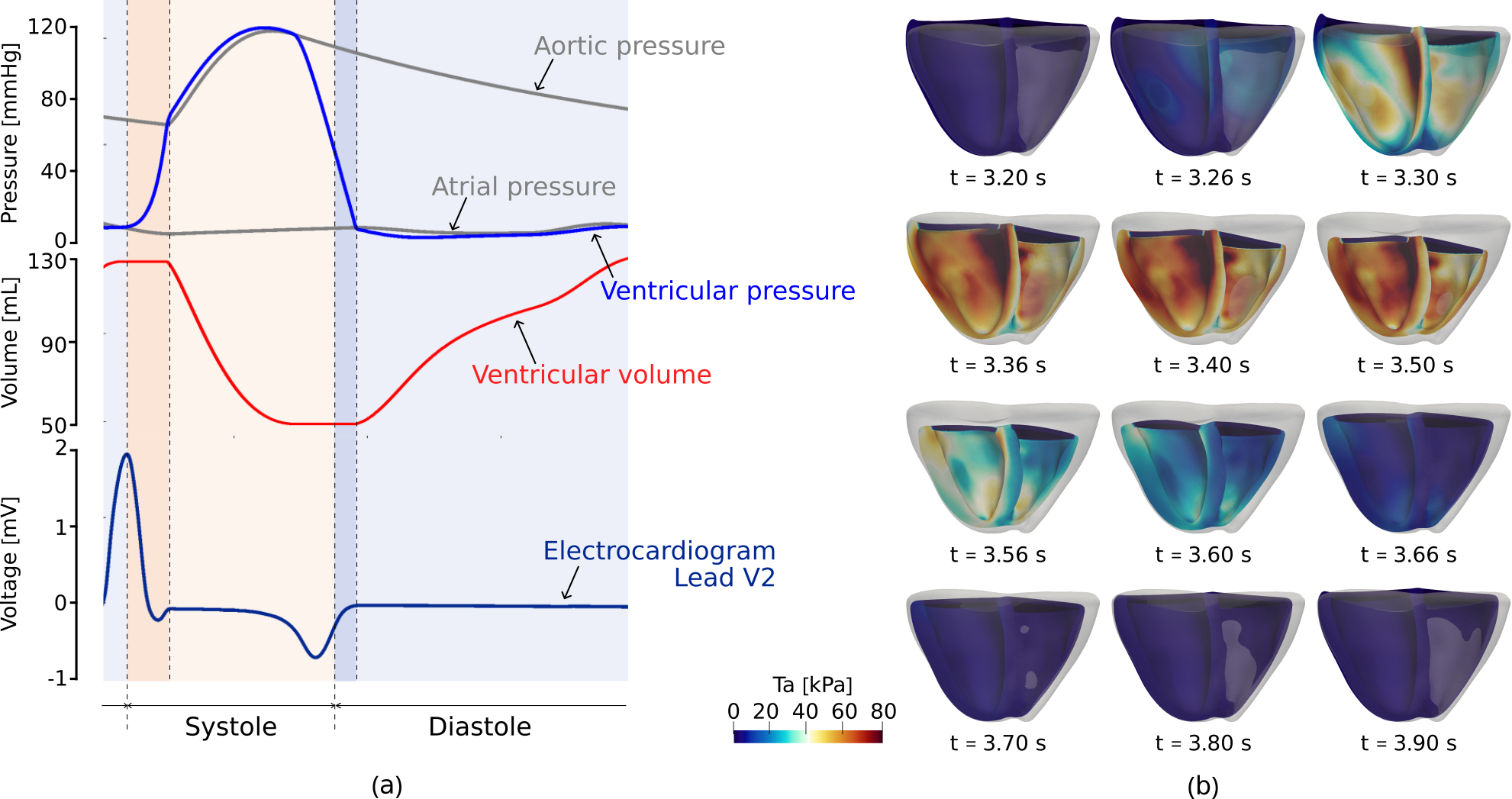}
	\caption{(a) Wiggers diagram obtained by solving the EMT problem when both heart and torso are in moving configuration. Background color identifies the four phases of the cardiac cycle: isovolumic contraction (dark orange), ejection (light orange), isovolumic relaxation (dark blue), and the remaining part of the diastole (light blue). (b) Cardiac mechanical deformation and active tension $T_a$ computed with the EMT model.}
	\label{Fig.Wigger}
\end{figure}

\begin{figure}[!t]
	\centering
	\includegraphics[width=0.73\textheight]{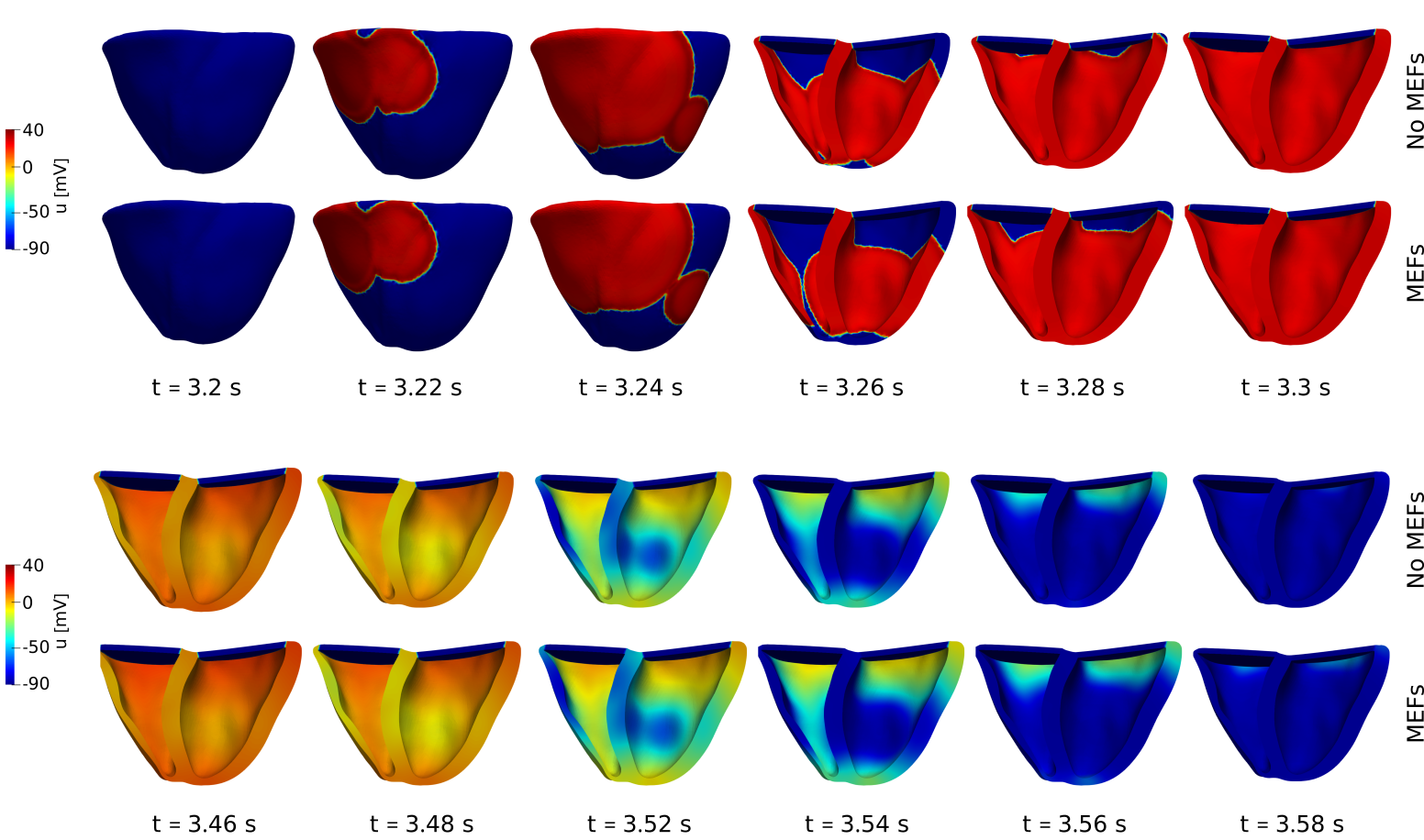}
	\caption{Propagation of the transmembrane potential $u$ on the cardiac reference domain $\Omega_H^0$ when the electrophysiological problem is solved on a static configuration (no MEFs, rows 1 and 3), and when is solved on a moving configuration (with MEFs, rows 2 and 4.)}
	\label{Fig.healthy_analysis_on_volume}
\end{figure}

Healthy conditions are achieved by using the model parameter setting reported in \cite{FEDELE2023115983}, specifically tailored for the biventricular geometry. The obtained results, encompassing quantitative indices related to cardiac mechanics and blood circulation serve as a proof of our model capability to simulate realistic mechanical contraction are reported in Figure \ref{Fig.Mechanics_results_loops}. Moreover, our EMT model allows for the integration of mechanical and electrophysiological aspects, yielding clinically relevant outputs and measurements. The Wiggers diagram in Figure \ref{Fig.Wigger} shows the alignment between systolic and diastolic phases, showcasing the heart contraction, ejection, and relaxation phases synchronized with the simulated ECG.

As illustrated in Figure \ref{Fig.healthy_analysis_on_volume}, we compare the transmembrane potential obtained on the deforming and static domain by projecting it on the static reference configuration. While the signal shape remains consistent, a noticeable reduction in the conduction velocities is observed when the MEFs are present. This is evident both in the depolarization and repolarization phases. 

\begin{figure}[!t]
	\centering
	\includegraphics[width=0.7\textheight]{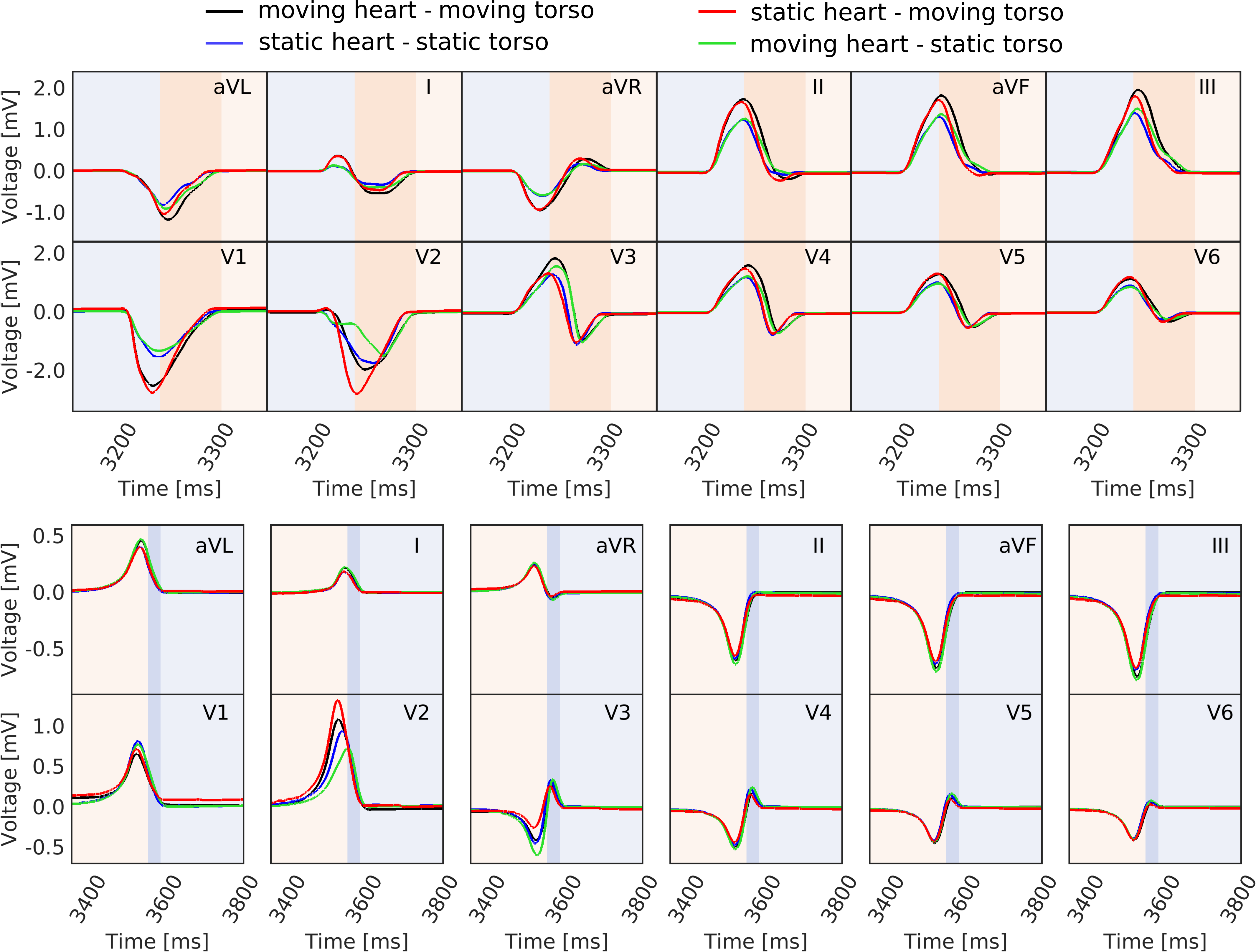}
	\caption{QRS (top) and T wave (bottom) computed with the EMT model with the four configurations described in Section \ref{Sub:comb_conf} and Table \ref{Tab.table_configurations_EMT}. Background color depicting four phases of the cardiac cycle: isovolumic contraction (dark orange), ejection (light orange), isovolumic relaxation (dark blue), and the remaining part of the diastole (light blue). }
	\label{Fig.ECG}
\end{figure}

\begin{figure}[!t]
	\centering
	\includegraphics[width=0.7\textheight]{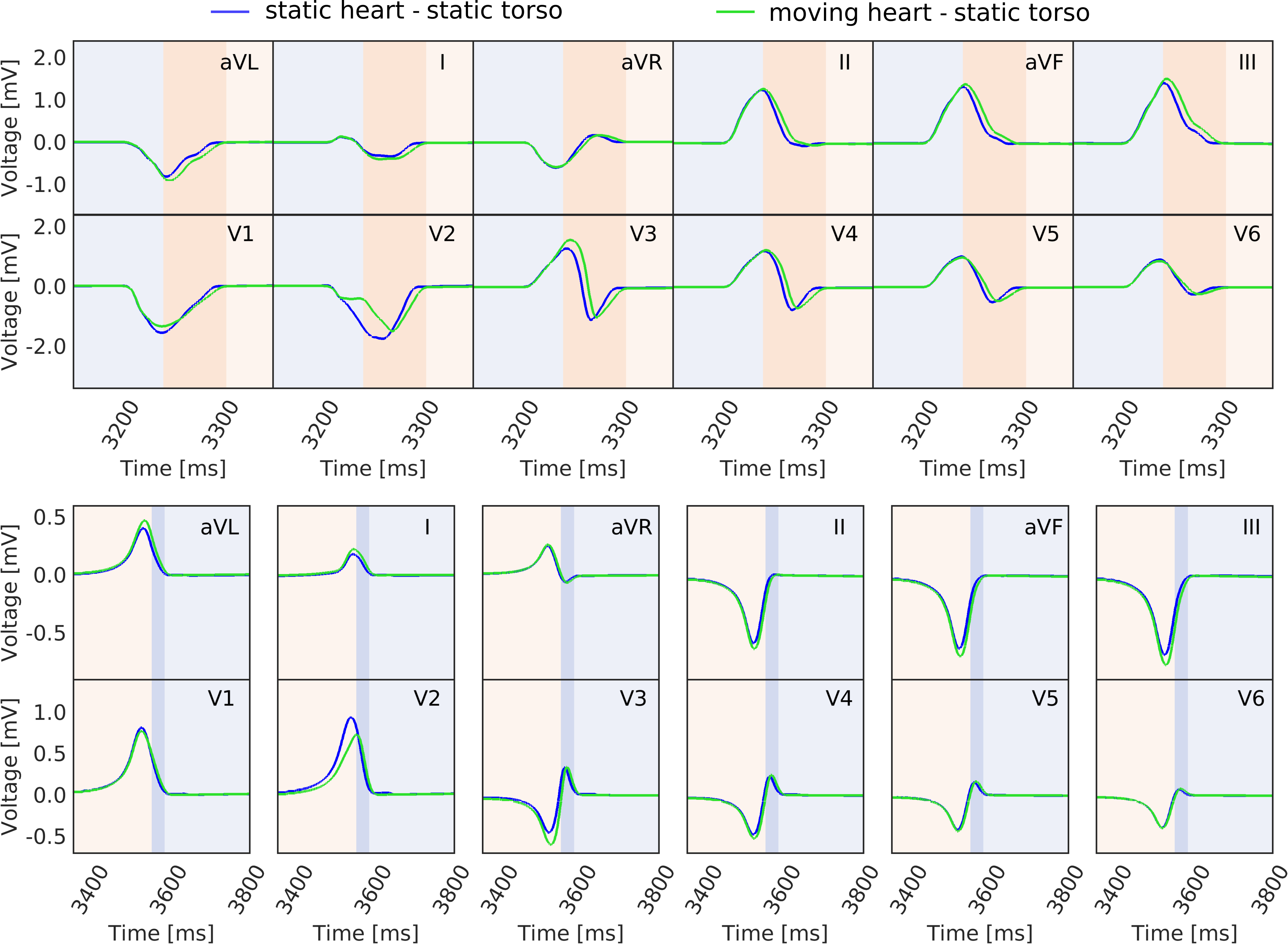}
	\caption{QRS (top) and T wave (bottom) computed with the EMT model when the torso is considered as a static domain, and the heart is either static or dynamic. Background color identifies four phases of the cardiac cycle: isovolumic contraction (dark orange), ejection (light orange), isovolumic relaxation (dark blue), and the remaining part of the diastole (light blue). }
	\label{Fig.ECG_healthy_static_torso}
\end{figure}

\begin{figure}[!t]
	\centering
	\includegraphics[width=0.7\textheight]{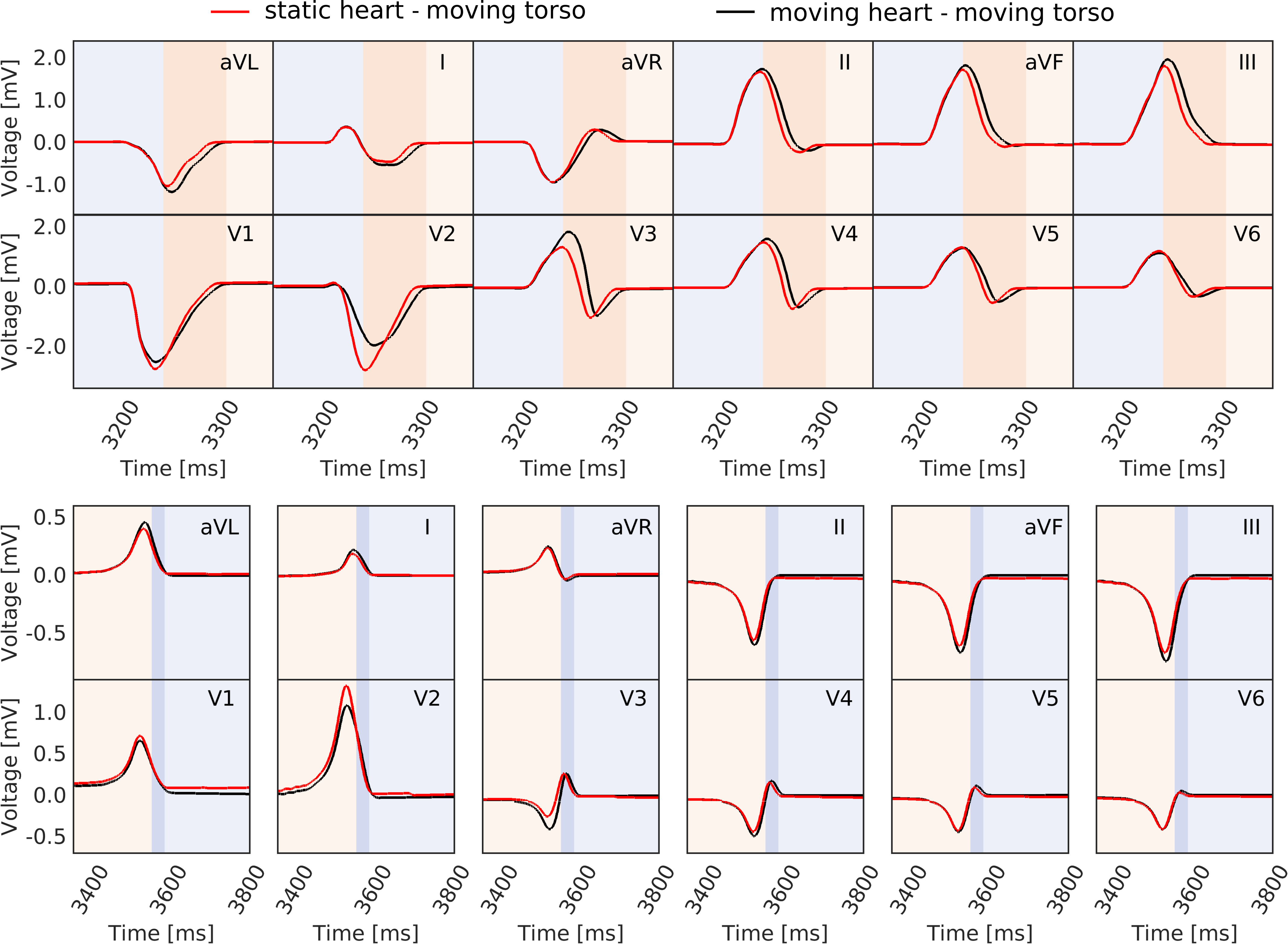}
	\caption{QRS (top) and T wave (bottom) computed with the EMT model when the torso is considered as a moving domain, and the heart is either static or dynamic. Background color identifies the four phases of the cardiac cycle: isovolumic contraction (dark orange), ejection (light orange), isovolumic relaxation (dark blue), and the remaining part of the diastole (light blue). }
	\label{Fig.ECG_healthy_moving_torso}
\end{figure}

In Figure \ref{Fig.ECG}, we depict a physiological QRS progression on the 12-lead ECG. Although the EP model lacks of heterogeneous apico-basal and transmural conduction velocities, which results in a non-physiological T wave, 
observed variations in the T wave due by myocardial contraction are independent from this heterogeneity, and can be therefore equally appreciated and investigated through this study.

Variations in transmembrane potential $u$ and torso potential $u_T$ are reflected in the corresponding ECG signals. In Figure \ref{Fig.ECG}, the QRS complex and T wave obtained from four different simulation settings are presented. Although the overall shape remains consistent for each ECG, differences in waves shape and shift are noticeable in all leads. 

Cardiac MEFs introduce a shift within the QRS wave in most leads, as observed when comparing signals obtained when the torso domain configuration is prescribed (see also Figure \ref{Fig.ECG_healthy_static_torso} and \ref{Fig.ECG_healthy_moving_torso}). Although this shift corresponds to a prolonged QRS duration, the ECG waves remain aligned with the cardiac contraction, ejection, and relaxation phases observed in the EM simulations.

\begin{figure}[!t]
	\centering
	\includegraphics[width=0.7\textheight]{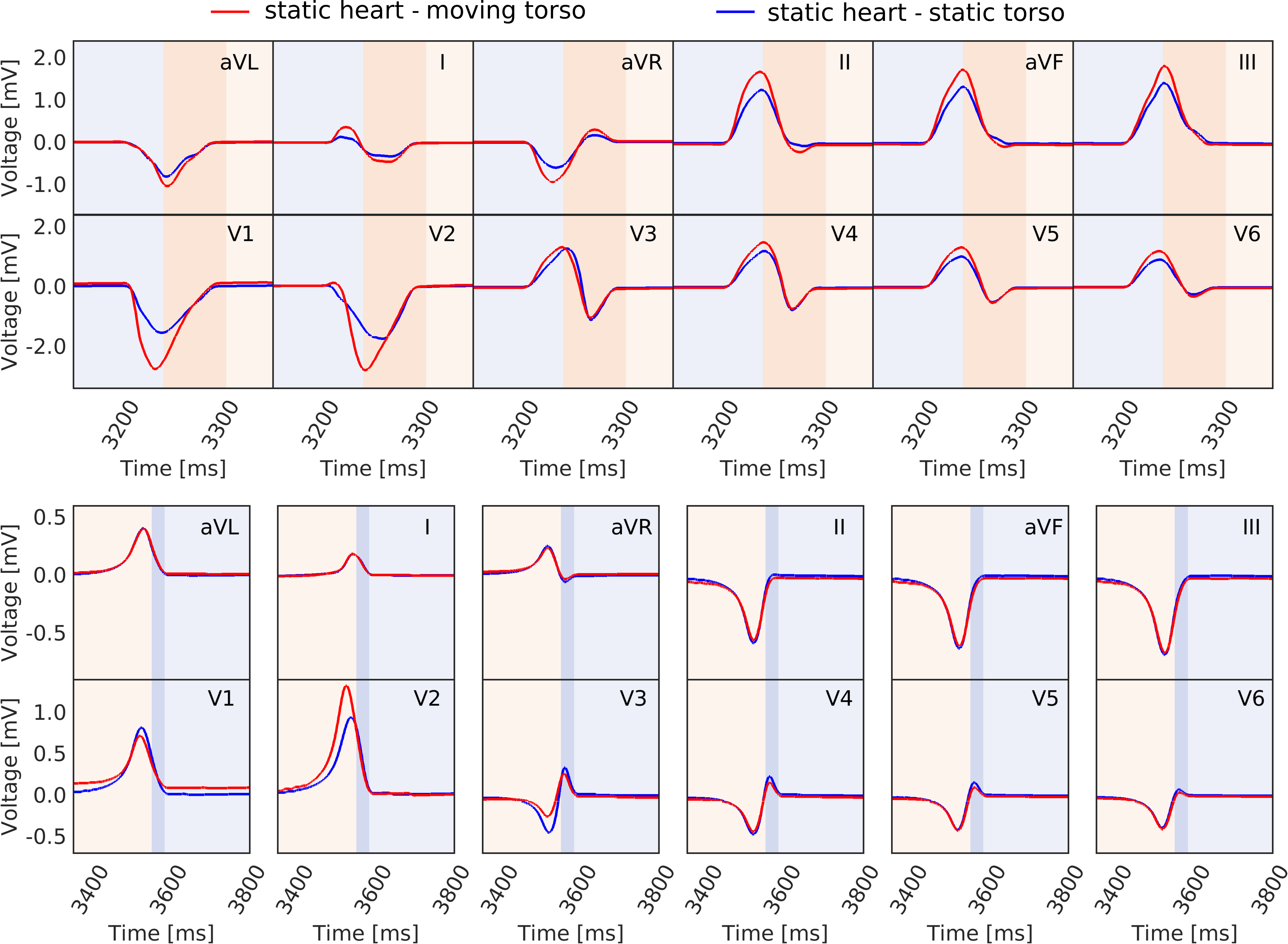}
	\caption{QRS (top) and T wave (bottom) computed with the EMT model when the heart is considered as a static domain, and the torso is either static or dynamic. Background color identifies the four phases of the cardiac cycle: isovolumic contraction (dark orange), ejection (light orange), isovolumic relaxation (dark blue), and the remaining part of the diastole (light blue). }
	\label{Fig.ECG_healthy_static_heart}
\end{figure}

\begin{figure}[!t]
	\centering
	\includegraphics[width=0.7\textheight]{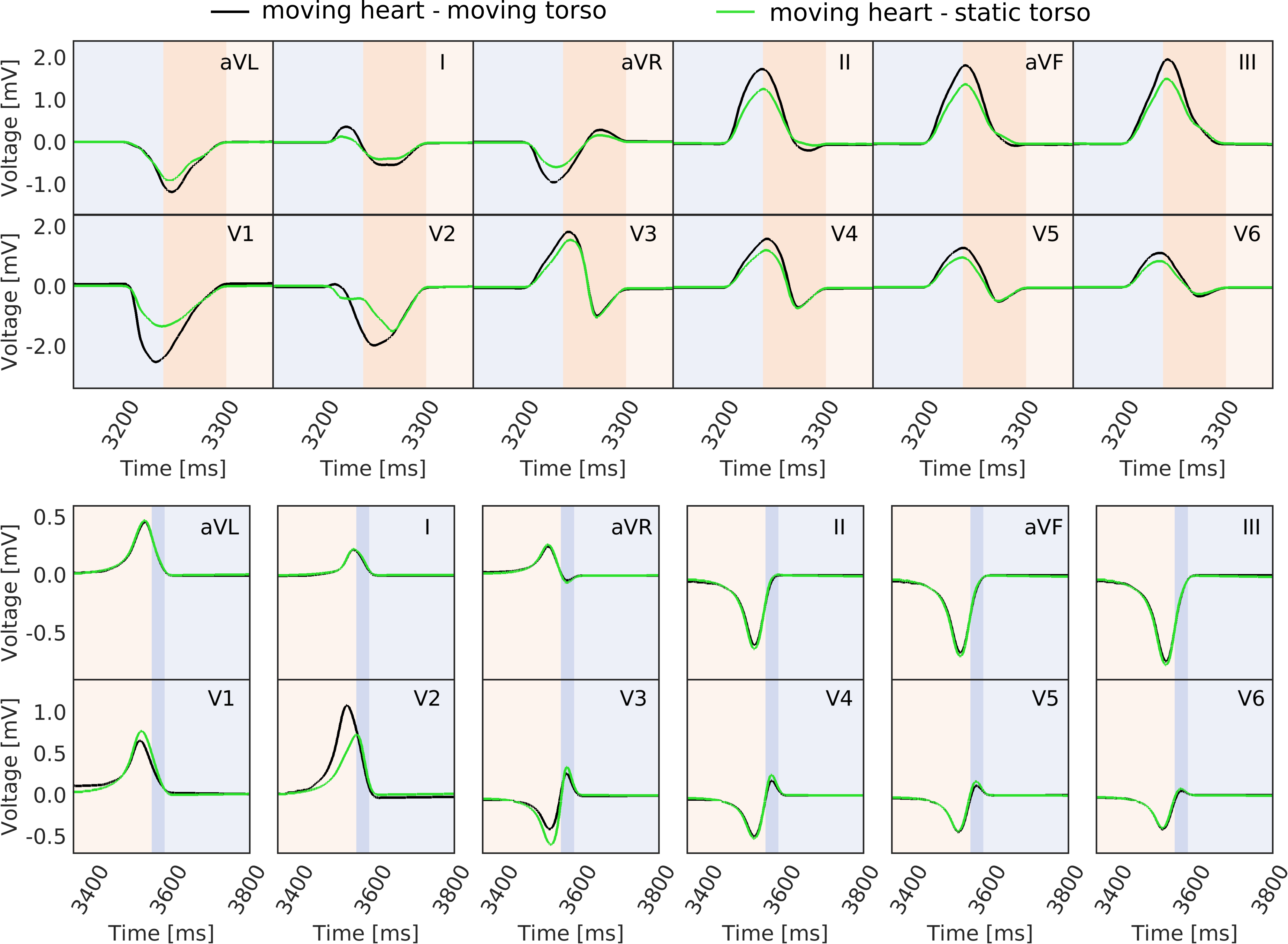}
	\caption{QRS (top) and T wave (bottom) computed with the EMT model when the heart is considered as a moving domain, and the torso is either static or dynamic. Background color identifies the four phases of the cardiac cycle: isovolumic contraction (dark orange), ejection (light orange), isovolumic relaxation (dark blue), and the remaining part of the diastole (light blue). }
	\label{Fig.ECG_healthy_moving_heart}
\end{figure}

Variations in torso potential, influenced by the deforming torso domain, result in changes in the amplitude of QRS waves. This effect is particularly prominent in limb leads $II$, $III$, $aVR$, and $aVF$. These leads are directly related to the left part of the torso, where the biventricular geometry is located, and consequently, the ECG and BSPMs are more influenced by the heart-torso interface deformation. Furthermore, substantial differences are observed in leads $V_{1}$, $V_2$ and $V_3$, which are closer to the cardiac domain and most susceptible to relative deformations (see also Figure \ref{Fig.ECG_healthy_static_heart} and \ref{Fig.ECG_healthy_moving_heart}).

The T wave, corresponding to the repolarization phase, is generally less affected by the four different simulation settings. However, leads $V_2$ and $V_3$ still exhibit changes in amplitude corresponding to the maximal contraction of the heart.

\begin{figure}[!t]
	\centering
	\includegraphics[width=0.55\textheight]{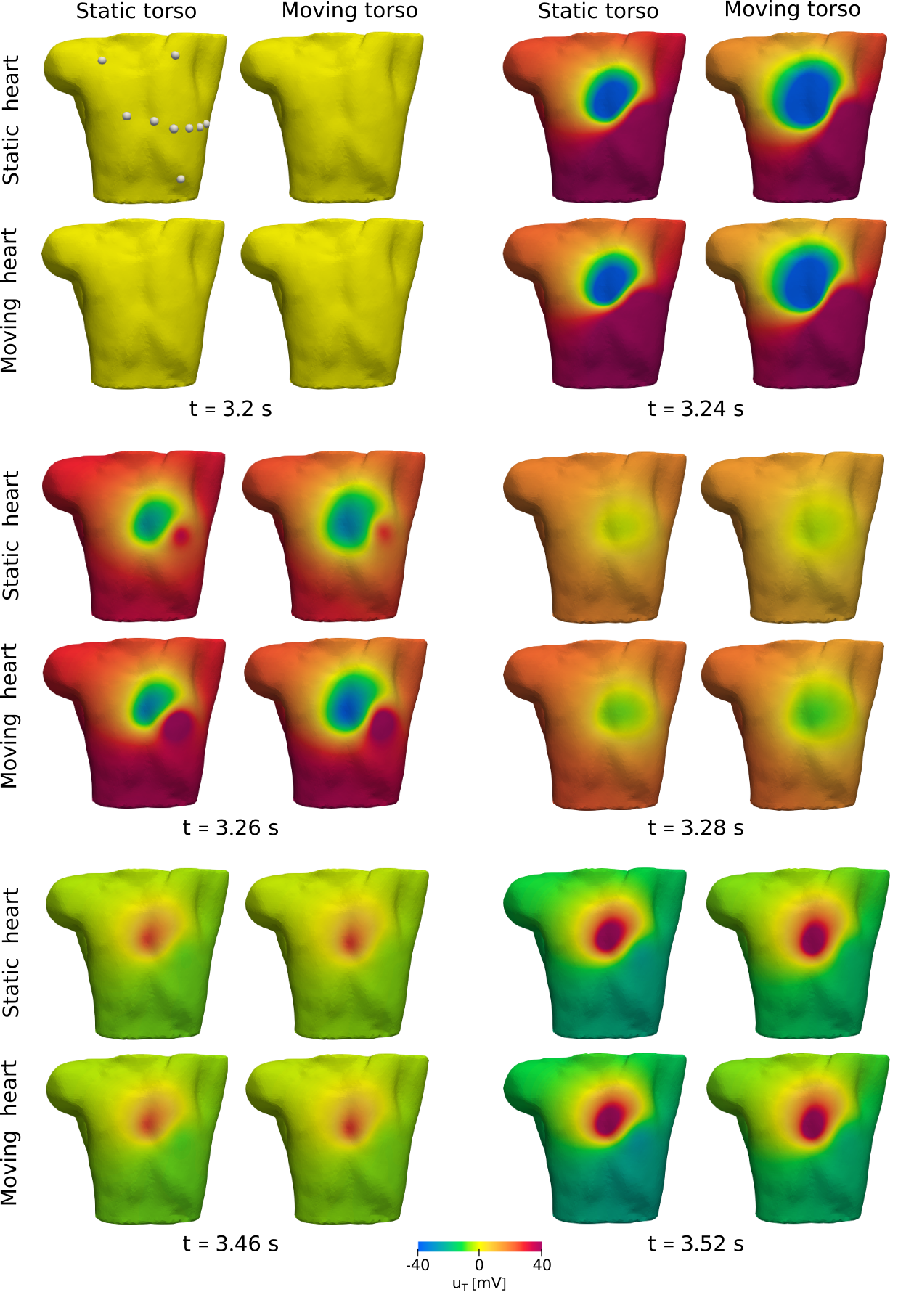}
	\caption{BSPMs on selected time instants for all four combinations of cardiac and torso configurations reported in Table \ref{Tab.table_configurations_EMT}. The electrodes locations are also reported.}
	\label{Fig.BSPM_healthy}
\end{figure}

Differences in signal propagation in the heart are also evident in the BSPMs depicted in Figure \ref{Fig.BSPM_healthy}. Starting from $t =$ \SI{3.26}{\second}, variations in BSPMs become significant when the torso configuration is prescribed, with fluctuations based on different configurations of the cardiac domain (as in Table \ref{Tab.table_configurations_EMT}, this is observed by looking at the group of four BSPMs by columns).
However, variations in BSPMs are not solely attributed to the cardiac MEFs. Different outcomes emerge when deforming the torso domain according to the cardiac displacement, instead of prescribing a static torso domain, as depicted in Figure \ref{Fig.BSPM_healthy} (as in Table \ref{Tab.table_configurations_EMT}, this is  observed by looking at the group of four BSPMs by rows).

To quantify differences in shape of two simulated ECG signals $\phi_1(i)$ and $\phi_2(i)$, $i = 1,\dots,N_T$, we use the following linear correlation coefficient (CC)
\begin{equation}
\label{Eq:cc}
\text{CC} = \frac{1}{s_1 s_2} \sum_{i=1}^{N_T}\left[\phi_1(i) - \overline{\phi_1}\right]\left[\phi_2(i) - \overline{\phi_2}\right],
\end{equation}
for each lead and its average across all leads. Here, $s_{1}$ and $s_{2}$ are the standard deviation of $\phi_{1}$ and $\phi_2$, respectively, while $\overline{\phi}_{1}$ and $\overline{\phi}_{2}$ are the corresponding arithmetic average values. The CCs between signals, obtained by reading Table \ref{Tab.table_configurations_EMT} by columns and rows, are reported in Tables \ref{Tab.cc_healthy_QRS} and \ref{Tab.cc_healthy_T}. This highlights a shift in the ECG wave, with the highest correlation coefficients corresponding to the heart-torso configuration with a prescribed torso.

\subsection{Pathological scenario}
\label{Sub:VT_conditions}
\begin{figure}[!t]
	\centering
	\includegraphics[width=0.6\textheight]{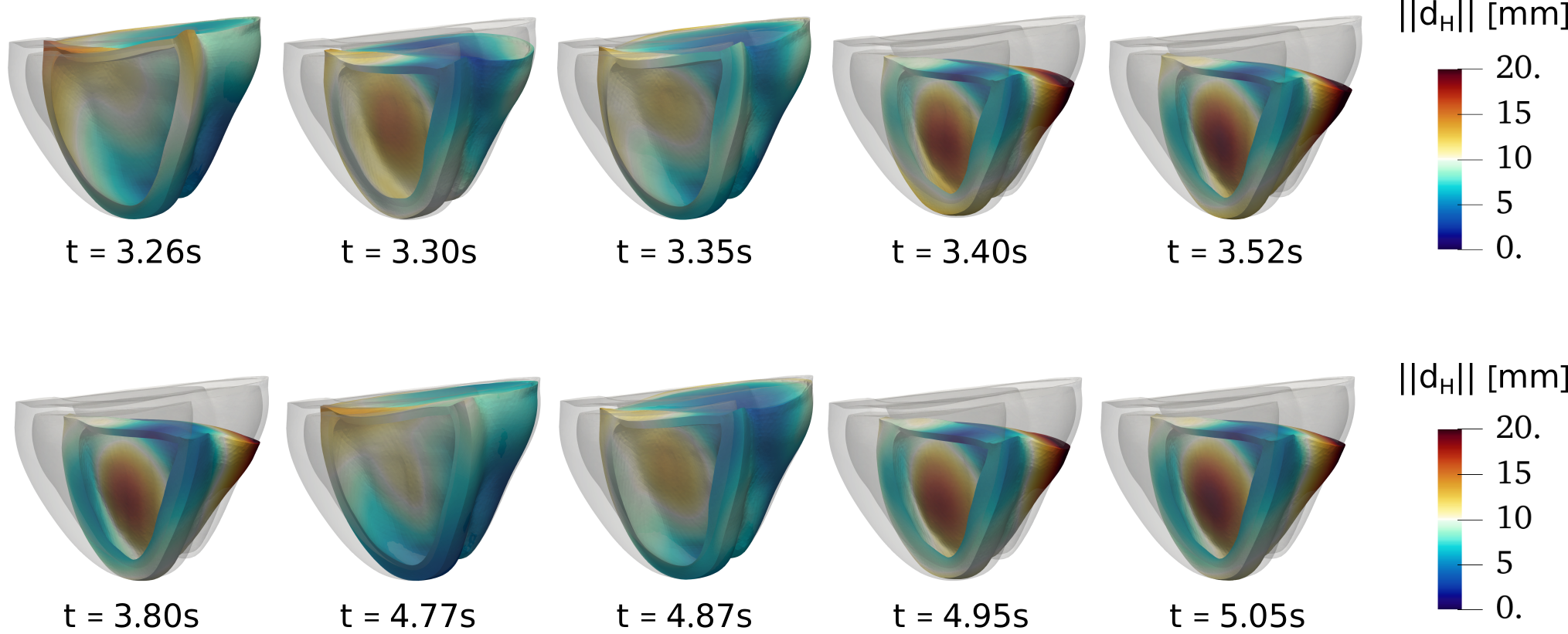}
	\caption{Displacement $\mathbf{d}_H$ simulated with the EMT model and MEFs on the heart.}
	\label{Fig.VT_displacement}
\end{figure}

Sustained VT is simulated in the biventricular geometry through an idealized septal ischemia. The VT is induced using an S1-S2-S3-S4 stimulation protocol, where different stimuli are applied at $t=$\SI{0}{\second}, $t=$\SI{0.45}{\second}, $t=$\SI{0.75}{\second} and $t =$ \SI{1.02}{\second} on the top part of the isthmus \cite{salvador2022role}.
To the best of our knowledge, this is the first time in which a VT is simulated using a fully-coupled 3D-0D EMT model in a biventricular geometry. In Figure \ref{Fig.VT_displacement}, we display different snapshots of displacement magnitude by considering the heart in the dynamic configuration. However, as we are primarily focused on observing electrophysiological bio-markers, we will limit our analysis to ECG and BSPMs, while employing the EMT model with different heart-torso configurations.

\begin{figure}[!t]
	\centering
	\includegraphics[width=0.6\textheight]{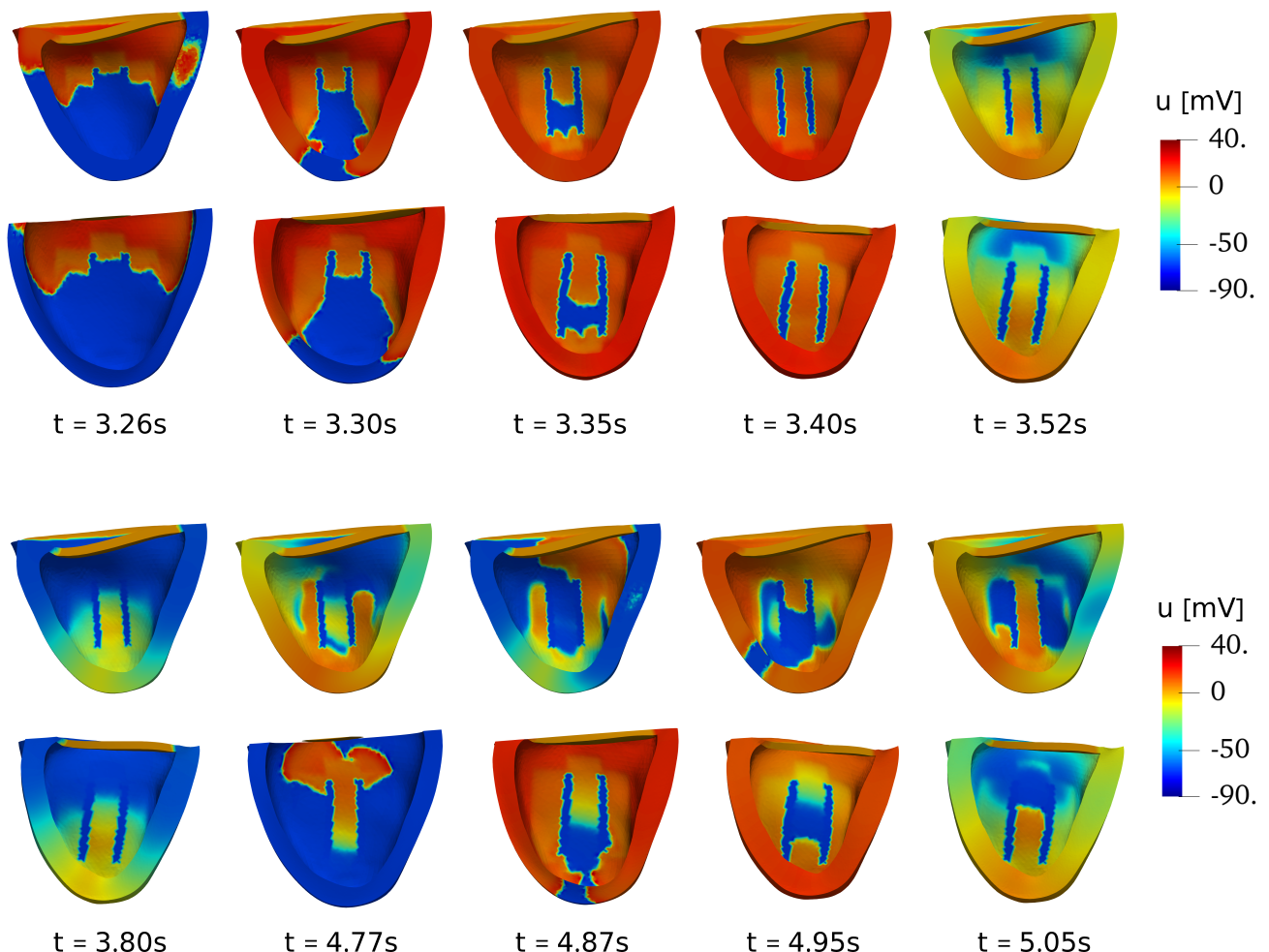}
	\caption{Propagation of the transmembrane potential $u$ on the heart $\Omega_H$ when the electrophysiological problem is solved on a static configuration (no MEFs, rows 1 and 3), and when is solved on a moving configuration (with MEFs, rows 2 and 4). When MEFs are activated, the cardiac geometry is warped by vector to show the deformation caused by the electro-mechanical interaction.}
	\label{Fig.signal_propagation_VT}
\end{figure}

\begin{figure}[!t]
	\centering
	\includegraphics[width=0.7\textheight]{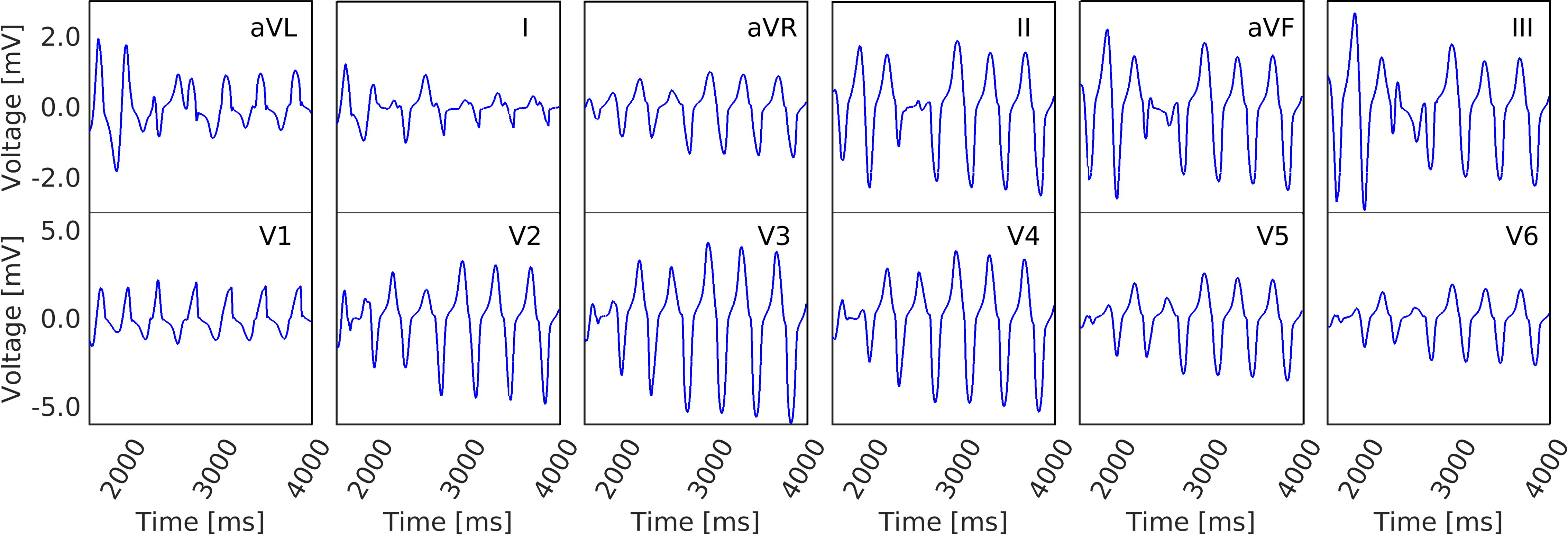}
	\caption{12-lead ECG computed considering both cardiac and torso domain in static configurations. The ECGs represent a polymorphic VT .}
	\label{Fig.ECG_VT_static_static}
\end{figure}

\begin{figure}[!t]
	\centering
	\includegraphics[width=0.7\textheight]{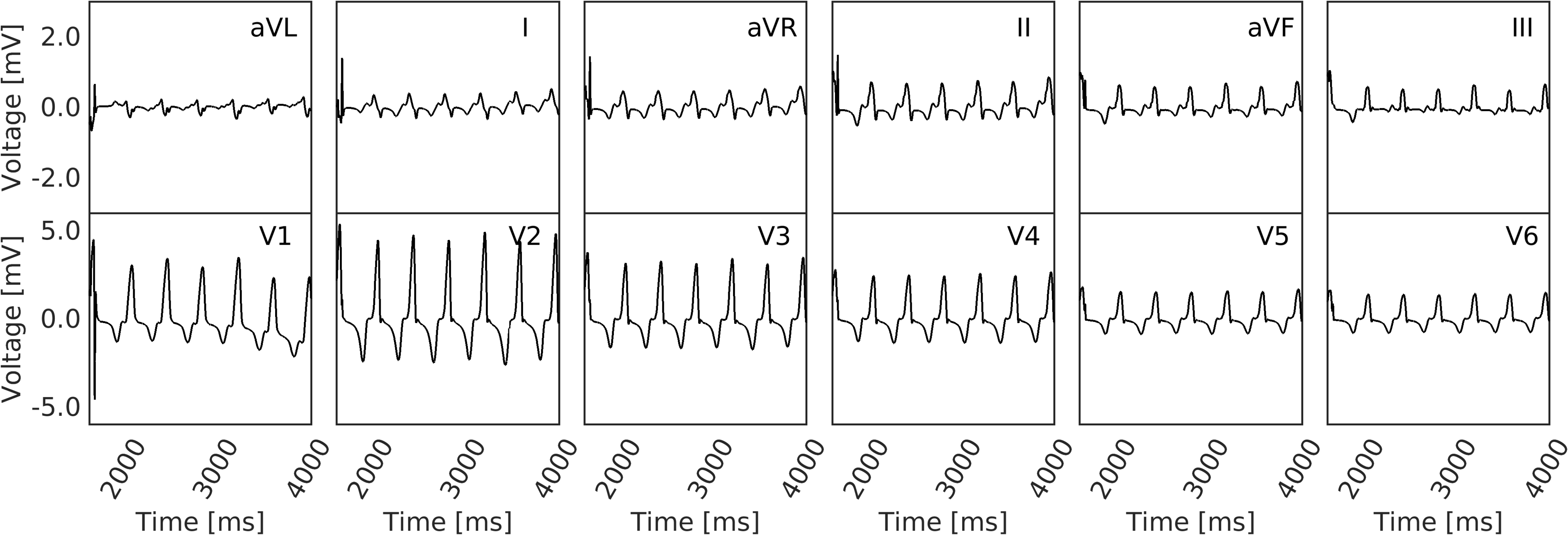}
	\caption{12-lead ECG computed considering both cardiac and torso domain in dynamic configurations. The ECGs represent a monomorphic VT .}
	\label{Fig.ECG_VT_moving_moving}
\end{figure}

In the EM simulation, we observe the formation and sustainment of a polymorphic VT. On the other hand, running an EMT simulation with a moving cardiac domain induces a stable monomorphic VT. This is evident in both the transmembrane potential propagation on the cardiac domain, as shown in Figure \ref{Fig.signal_propagation_VT}, and the ECG traces presented in Figures \ref{Fig.ECG_VT_static_static} and \ref{Fig.ECG_VT_moving_moving}.


We compare BSPMs and ECGs obtained by solving the EMT model with the heart in a dynamic domain, and the torso in either static or dynamic configuration. 

\begin{figure}[!t]
	\centering
	\includegraphics[width=0.6\textheight]{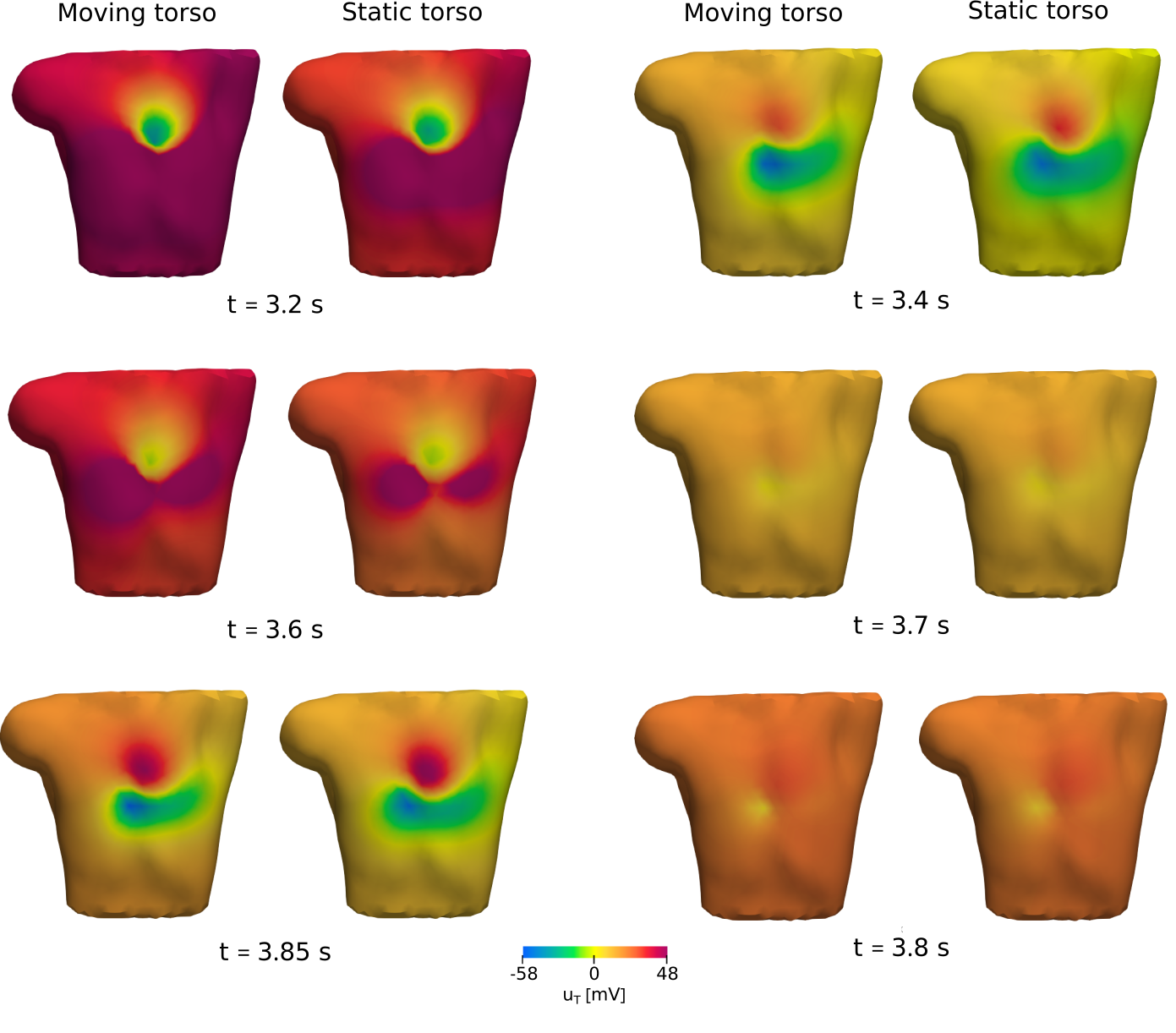}
	\caption{BSPMs on selected time instants computed with a prescribed cardiac moving configuration according to the second row of Table \ref{Tab.table_configurations_EMT}. 	}
	\label{Fig.VT_BSPM}
\end{figure}

BSPMs for the static and moving configurations of the heart are illustrated in Figure \ref{Fig.VT_BSPM}. Differences in signal magnitude are present, while the overall signal distribution remains consistent across the entire torso surface.

\begin{figure}[!t]
	\centering
	\includegraphics[width=0.7\textheight]{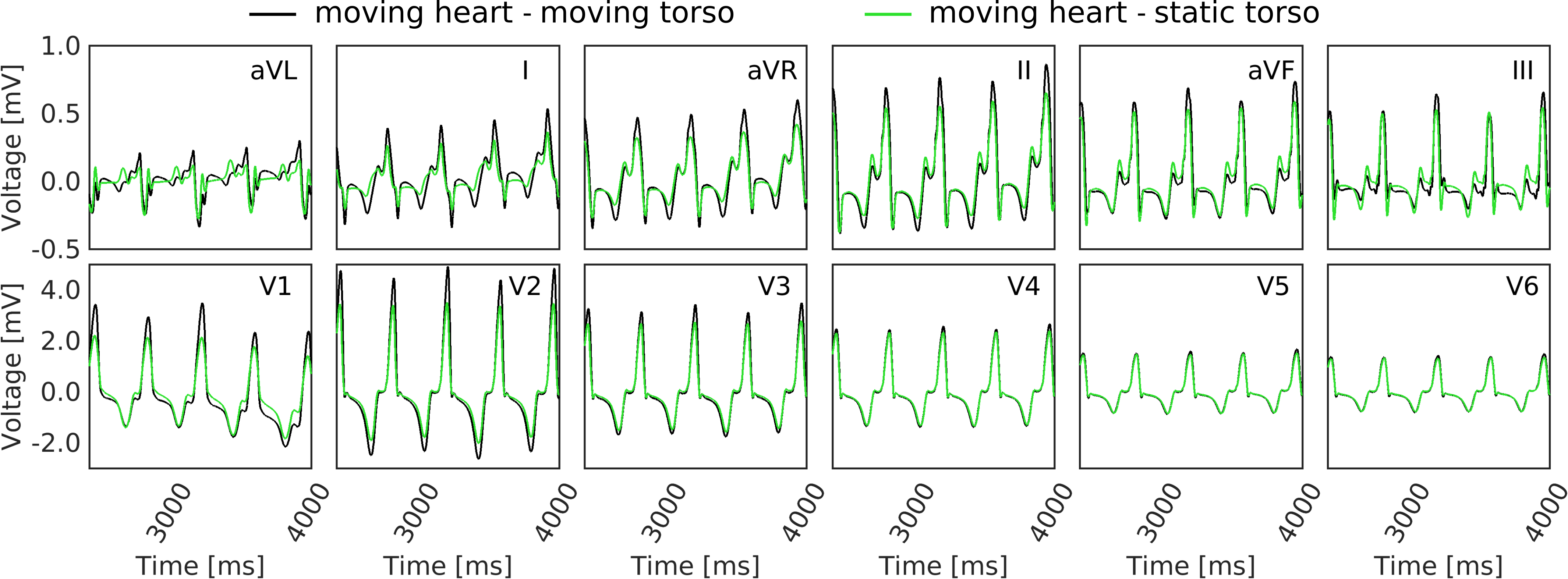}
	\caption{Comparison of 12 lead ECGs when the EMT model is solved in a cardiac moving configurations. The ECGs represent a monomorphic VT .}
	\label{Fig.ECG_VT_moving_heart}
\end{figure}

The ECG traces, depicted in Figure \ref{Fig.ECG_VT_moving_heart}, showcase variations in the amplitude of both the QRS complex and T-wave in all limb leads and leads $V_1$, $V_2$, $V_3$, with other precordial leads remaining largely unchanged.   Alteration in T-wave shape are mostly noticeable in the limb leads, particularly in lead $I$, where T-wave inversion is consistently observed throughout all heartbeats.

The CCs over the entire set of leads are reported in Table \ref{Tab.cc_VT_QRS}, confirming the substantial impact of torso domain pseudo-deformation on the ECG, especially in leads $I$ and $II$. The average CC is also generally smaller than the corresponding values obtained in healthy simulations, showing a greater shape variation in pathological simulations rather than healthy ones.

\begin{landscape}
	\begin{table}
		\centering
		\begin{tabular}{c|cccccccccccc|c}
			\hline
			&\multicolumn{12}{c}{CC}&\\\hline
			Heart-torso vs heart-torso &aVL &I & aVR &II &aVF &III &V$_1$ &V$_2$ &V$_3$ &V$_4$ &V$_5$ &V$_6$ &mean \\ \hline
			Static-static vs moving-static  &0.9763 &0.9763 &0.9847 &0.9811 &0.9786 &0.9770 &0.9911 &\cellcolor{red!20}0.8550 &0.8765 &0.9425 &0.9677 &0.9864 &0.9578 \\
			Moving-moving vs static-moving &0.9723 &0.9779 &0.9791 &0.9746 &0.9723 &0.9711 &0.9939 &0.9209 &\cellcolor{red!20}0.8216 &0.9316 &0.9616 &0.9817 &\cellcolor{red!20}0.9549\\
			Moving-moving vs moving-static &0.9942 &0.9624 &0.9981 &0.9973 &0.9978 &0.9979 &0.9758 &0.8682 &\cellcolor{green!25}0.9926 &0.9907 &0.9924 &0.9988 &0.9805\\
			Static-static vs static-moving &0.9946 &0.9629 &0.9950 &0.9960 &0.9974 &0.9982 &0.9703 &\cellcolor{green!25}0.9360 &0.9668 &0.9937 &0.9958 &0.9996 &\cellcolor{green!25}0.9839\\
			 \hline		
		\end{tabular}
		\caption{QRS correlation coefficients (CC) in the standard 12 leads of a healthy patient between different heart-torso configurations. The first two rows correspond to a column-wise interpretation of Table \ref{Tab.table_configurations_EMT}, while rows 3 and 4 represent a row-wise reading  of Table \ref{Tab.table_configurations_EMT}. The average CC across all leads is included in the last column. For leads with largest ECG shape variations, we highlighted in red the configuration with smallest CC (highest variation) and in green the configuration with largest CC (smaller variation).}
		\label{Tab.cc_healthy_QRS}
\end{table}
\begin{table}
	\centering
	\begin{tabular}{c|cccccccccccc|c}
		\hline
		&\multicolumn{12}{c}{CC}&\\\hline
		Heart-torso vs heart-torso &aVL &I & aVR &II &aVF &III &V$_1$ &V$_2$ &V$_3$ &V$_4$ &V$_5$ &V$_6$ &mean \\ \hline
		Static-static vs moving-static &0.9935 &0.9939 &0.9952 &0.9943 &0.9939 &0.9936 &0.9968 &0.9572 &0.9517 &0.9839 &0.9912 &0.9966 &0.9868\\
		Moving-moving vs static-moving &0.9928 &0.9921 &0.9909 &0.9919 &0.9923 &0.9926 &0.9886 &\cellcolor{green!25}0.9905 &\cellcolor{red!20}0.8598 &0.9772 &0.9894 &0.9954 &\cellcolor{red!20}0.9795\\
		Moving-moving vs moving-static &0.9985 &0.9982 &0.9955 &0.9968 &0.9975 &0.9979 &0.9808 &\cellcolor{red!20}0.9198 &\cellcolor{green!25}0.9853 &0.9959 &0.9916 &0.9941 &0.9877\\
		Static-static vs static-moving &0.9993 &0.9973 &0.9991 &0.9993 &0.9993 &0.9993 &0.9934 &0.9614 &0.9263 &0.9988 &0.9952 &0.9966 &\cellcolor{green!25}0.9888\\
		 \hline		
	\end{tabular}
	\caption{T wave correlation coefficients (CC) in the standard 12 leads of a healthy case between different heart-torso configurations. The first two rows correspond to a column-wise interpretation of Table \ref{Tab.table_configurations_EMT}, while rows 3 and 4 represent a row-wise reading  of Table \ref{Tab.table_configurations_EMT}. The average CC across all leads is included in the last column. For leads with largest ECG shape variations, we highlighted in red the configuration with smallest CC (highest variation) and in green the configuration with largest CC (smaller variation).}
	\label{Tab.cc_healthy_T}
\end{table}
\begin{table}
	\centering
	\begin{tabular}{c|cccccccccccc|c}
		\hline
		&\multicolumn{12}{c}{CC}&\\\hline
		Heart-torso vs heart-torso &aVL &I & aVR &II &aVF &III &V$_1$ &V$_2$ &V$_3$ &V$_4$ &V$_5$ &V$_6$ &mean \\ \hline
		Moving-moving vs static-moving &\cellcolor{red!20}0.7763 &\cellcolor{red!20}0.9180 &0.9676 &0.9658 &0.9506 &0.9076 &0.9807 &0.9840 &0.9915 &0.9985 &0.9986 &0.9986 &0.9531\\ \hline		
	\end{tabular}
	\caption{Correlation coefficients (CC) between the standard 12 leads of a pathological case, computed when the heart is considered in a moving domain. The average CC across all leads is included in the last column. Largest variations are highlighted with red color.}
	\label{Tab.cc_VT_QRS}
\end{table}
\vspace{4cm}
\end{landscape}

	\section{Discussion}
\label{Sec:discussion}
In the simulation of cardiac healthy conditions, our model demonstrates its ability to accurately represent physiological scenarios in agreement with reference values from the medical literature of pressure and volume loops, ECGs and BSPMs. The interplay of mechanical deformation on the ECG and BSPMs is highlighted, revealing a direct link between cardiac MEFs and variations in ECG and BSPM signals. Differently than in \cite{moss2021fully}, our model features shifts and prolongations of QRS waves across all leads during the contraction and ejection phases. Additionally, amplitude variations are observed in ECG signals, particularly in the limb leads associated with the left torso and precordial leads $V_1$, $V_2$, and $V_3$, primarily attributed to torso domain deformation. T-wave variations are predominantly registered in lead $V_2$, in line with \cite{moss2021fully}, emphasizing changes in wave amplitude.

Regarding the pathological conditions, a VT is simulated by defining a sustained figure-of-eight reentry with an an idealized isthmus located at the (bi)ventricular septum. Moreover, we highlight the impact of the cardiac MEFs on the propagation of the electrical signal. Indeed, MEFs may transform the nature of the simulated VT from polymorphic to monomorphic, therefore altering ECG and BSPM patterns. While in this case cardiac MEFs have a prominent role in changing the EP signal, effect of the torso domain deformation has also been observed. Specifically, changes in amplitude and shape of the ECGs are shown, with major variations in lead $I$, where inverted polarity in the T wave is achieved.

Compare to the more recent EM computational model presented in \cite{moss2021fully},  our methodology still ensures complete equivalence in terms of physical modeling, as we address the same underlying mathematical problem, but features some differences in terms of computational methods and capacity of simulating pathological scenario. Specifically, in \cite{moss2021fully} the authors constructed a halo around the four-chamber heart within the torso geometry. For each new heart configuration defined by the solution of an electromechanical model \cite{gerach2021electro}, a new mesh on the halo was created to account for the altered heart shape. Subsequently, the BEM was solved within the halo-torso geometry to calculate the ECGs and the BSPMs. This modeling approach, therefore, necessitated a remeshing stage for each time instance of EP output computation. Moreover, due to the use of BEM, it lacked the capability to represent anisotropic conduction properties in the extracellular space, which became significant in pathological scenarios. Conversely, our approach avoids the necessity for remeshing, and allows for ECGs and BSPMs computation even for pathological conditions.

Furthermore, in the present work the influence of MEFs on EM outcomes results to be important and more pronounced than in the test cases presented in \cite{salvador2021electromechanical,REGAZZONI2022111083}. 
It therefore remains unclear how the chosen MEFs and cardiac geometry can impact the simulation outputs. Nevertheless, our findings indicate that both MEFs and cardiac geometry can lead to distinct outcomes, particularly in the propagation of EP signals. Further insights could be gained by testing the proposed EMT model on a cohort of four-chamber hearts, similar to the one presented in~\cite{strocchi2020publicly}. Moreover, numerous studies have demonstrated that the position, dimensions, and shape of the heart~\cite{NGUYEN2015617,Minchole2019}, along with the presence and configuration of organs within the torso \cite{Keller2010,Sanchez2018}, can impact the shape and amplitude of ECGs. This further highlight the need to expand the range of cardiac and torso domains considered in this analysis. 

Finally, the primary limitations of our approach reside in the representation of heterogeneous cardiac conduction properties, which influence T-wave shape and polarity, as well as a realistic representation of the Purkinje network. These two features will be implemented in the future, and their effect on the EMT simulation will be further investigated.

\section{Conclusions}
\label{Sec:conclusions}
In this study, we introduced an EMT model by unidirectionally coupling an EM model of the heart with a passive conductive model of the torso. This mathematical model dynamically defines the torso domain deformation resulting from myocardial displacement. The flexible segregated-intergrid-staggered numerical framework employed allows for the arbitrary and independent selection of the cardiac displacement used to solve both the Monodomain model and to modify the torso domain, which is equivalent to employ static or dynamic cardiac and torso domain configurations. This allows to explore the impact of myocardial deformation on EP propagation and, consequently, on the ECG and BSPMs, as well as investigate the influence of shifting the heart-torso surface, and thereby the torso domain, on EP outputs. 

The model is tested under both healthy and pathological scenarios, the latter involving cardiac arrhythmias, specifically VT. To ensure fairness in our comparison, as static configuration, the displacement of the cardiac geometry in the end-diastolic phase is extracted from the EM simulation and imposed to the Monodomain and torso lifting problem. Overall, based on the results obtained from our model, we concluded that the influence of cardiac contraction on EP outputs should not be underestimated, particularly when simulating pathological conditions. 
	\section*{Acknowledgements}
The present research is part of the activities of ``Dipartimento di Eccellenza 2023--2027'', MUR, Italy, Dipartimento di Matematica, Politecnico di Milano. R. Piersanti and L. Dede' have received support from the project PRIN2022, MUR, Italy, 2023--2025, 202232A8AN ``Computational modeling of the heart: from efficient numerical solvers to cardiac digital twins''. F. Regazzoni has received support from the project PRIN2022, MUR, Italy, 2023--2025, P2022N5ZNP ``SIDDMs: shape-informed data-driven models for parametrized PDEs, with application to computational cardiology''. The authors acknowledge their membership to INdAM GNCS - Gruppo Nazionale per il Calcolo Scientifico (National Group for Scientific Computing, Italy). This project has been partially supported by the INdAM-GNCS Project CUP E53C22001930001.


	\bibliography{references}
	
	\newpage
	\begin{appendix}
	\section{Cardiac electro-mechanical model}
\label{App:cardiac_EM}
In this Appendix we report the mathematical models employed  to define cardiac electromechanics and cardiovascular hemodynamics. These models are also briefly described and accounted for in Section \ref{Sec:models}. Finally, the 12-lead ECG system is described in the last section.

\subsection{Electrophysiological model}
The cardiac transmembrane potential $u$ is computed by solving the Monodomain model endowed with gradient deformation tensors to account for the the mechanical displacement on the potential propagation. Specifically, the formulation of the Monodomain model reads:
\begin{subequations}\label{eqn: EP}
	\begin{empheq}[left={\empheqlbrace\,}]{align} \begin{split}&\chi \left[ \EPCm \dfrac{\partial u}{\partial t} + I_{ion}(u,\mathbf{w},\mathbf{c}) \right]- \nabla \cdot ( J_H \mecF_H^{-1} \mathbf{D}_m \, \mecF_H^{-T} \nabla u) \\ &\qquad \qquad \:\, = J_H \EPchim I_{app}(t)  \qquad  \:\:\:\:\:\:\:\:\:\:\, \qquad\qquad \text{in } \{\Omega_{H}^0 \cup \Omega_C^0\} \times (0,T],\end{split} \label{eqn: EP1}\\
	& \dfrac{\partial \mathbf{w}}{\partial t} - \EPrhsGating(u, \mathbf{w}) = \boldsymbol{0}  \:\:\,  \qquad \qquad \qquad  \qquad
	\quad \:\:\:\:\:\:\:\, \text{in } \{\Omega_{H}^0\cup \Omega_C^0\} \times (0,T], \label{eqn: EP3}\\
	&\frac{d \mathbf{z}}{d t} = \mathbf{G}(u,\mathbf{w},\mathbf{z})  \qquad \qquad \qquad \qquad \qquad \quad  \quad \:\:\:\,\text{in } \{\Omega_H^0 \cup \Omega_C^0 \}\times (0,T],\label{eqn: EP5}	\\ 
	&(J_H \mecF_H^{-1}\mathbf{D}_m \mecF_H^{-T}\nabla u)\cdot \mathbf{n}_H = 0 \qquad \qquad  \quad \:\:\:\:\:\, \text{on } \partial \{\Omega_{H}^0 \cup \Omega_C^0\} \times (0, T], \label{eqn: EP4}\\
	&u = u_{m,0}, ~ \mathbf{w} = \mathbf{w}_{0}, ~ \mathbf{z}= \mathbf{z}_0 \qquad \qquad \qquad \:\:\:\:\:\:\,\text{in }\{\Omega_{H}^0 \cup \Omega_C^0\} \times \{t = 0\}.
	\end{empheq}
\end{subequations}
where $C_m$ is the capacitance per unit area, and $\chi$ is the surface-to-volume ratio of the membrane.
Here $I_{app}$ is a function representing the activation sites used as pacing protocol (referred to Section \ref{Sec:numerical_results}), and $\mathbf{D}_m$ is the conductivity tensor, computed as:
$$ \boldsymbol{D}_{m} = \mu \sigma_{\ell}^{m} \frac{\mecF_H\fZero \otimes \mecF_H\fZero}{\|\mecF_H\fZero\|^2} + \mu \sigma_{t}^{m} \frac{\mecF_H\sZero \otimes \mecF_H\sZero}{\|\mecF_H\sZero\|^2} + \mu \sigma_{n}^{m} \frac{\mecF_H\nZero \otimes \mecF_H\nZero}{\|\mecF_H\nZero\|^2},$$
The conduction coefficients $\sigma_{\ell,t,n}^m$ are strictly related to the extracellular and intracellular conduction coefficients by means of the following formula:
$$ \sigma_{\ell,t,n}^m = \frac{\sigma_{\ell,t,n}^i \sigma_{\ell,t,n}^e}{\sigma_{\ell,t,n}^i + \sigma_{\ell,t,n}^e}.$$
Equations \eqref{eqn: EP3}-\eqref{eqn: EP5} stand for the ionic model.

\subsection{Cardiomyocyte active contraction}
The cardiomyocyte active contraction is modelled with either the RDQ18 and RDQ20 models, which have the following structure:
\begin{subequations}
	\begin{empheq}[left={\empheqlbrace\,}]{align}
	&\frac{\partial \ActStateHF}{\partial t} = \ActRhs(\ActStateHF,\Cai,\SL,\frac{d\SL}{dt})  \qquad \:\:\:\:\:\:\ \qquad \text{in } \{\Omega_{H}^0\cup \Omega_C^0\} \times (0,T], \label{eqn: Act}\\
	&\mathbf{s} = \mathbf{s}_0 \qquad  \qquad \qquad \qquad \qquad \qquad \qquad \:\:\: \text{on }\{\Omega_{H}^0 \cup \Omega_C^0\}\times \{t = 0\},
	\end{empheq}
\end{subequations}
where the unknown is the vector $\mathbf{s}$ of state variables, $\mathbf{K}$ a suitable function (see \cite{regazzoni2018redODE}), and $SL$ is obtained from the mechanical model as:
$$ SL = SL_0\sqrt{\mathcal{I}_{4f}(\mathbf{d}_H)}.$$
The variable $SL_0$ is the sarcomeres length at rest, while $\mathcal{I}_{4f} = \mathbf{F}_H\mathbf{f}_0 \cdot \mathbf{F}_H\mathbf{f}_0$ is a measure of the tissue stretch along the fibers direction. 

The active tension $T_a$ is computed as:
$$ T_a(\mathbf{s}) = T_a^{\mathrm{max}}G(\mathbf{s},P)\left[\hat{\xi} + C_{LRV}(1-\hat{\xi}) \right],$$
where $P$ is the permissivity, that is the  fraction of contractile units being in the force-generation state, 
$T_a^{\mathrm{max}}$ is the total tension generated (obtained when $P=1$), $G(\mathbf{s})$ is a linear function related to the permissivity~ \cite{regazzoni2020biophysically}, $\hat{\xi}\in[0,1]$ is the normalized intra-ventricular distance, and $C_{LRV}\in (0,1]$ is the left-right ventricular contractility ratio.

\subsection{Cardiac active and passive contraction}
The large deformation of the cardiac tissue is described through the finite elasticity theory. 
The displacement $\mathbf{d}_H$ is obtained by solving the momentum conservation equation endowed with proper boundary conditions and a formulation for of the Piola-Kirchhoff stress tensor $\mathbf{P}$. This yields to the following model:
\begin{subequations}	\begin{empheq}[left={\empheqlbrace\,}]{align}
	& \rho_{\text{s}} \dfrac{\partial^2 \displ_H}{\partial t^2} - \nabla \cdot \tenspiola(\displ_H \Tens(\ActStateHF)) = \boldsymbol{0}
	\qquad \qquad \:\:\:\:\:\:\:\ \text{in } \{\Omega_{H}^0 \cup \Omega_{C}^0\} \times (0,T], \label{eqn: Mec1}\\
	& \tenspiola(\displ_H, \Tens(\ActStateHF)) \mecNref = \BCmecKepiTens \displ_H +
	\BCmecCepiTens \dfrac{\partial \displ_H}{\partial t}
	\qquad \qquad \:\:\:\:\:\:\:\ 
	\text{on } \Gamma_H^{\text{epi}} \times (0, T], \label{eqn: Mec2}\\
	&\tenspiola(\displ_H, \Tens(\ActStateHF)) \mecNref = -\PLV(t) \, J \mecF^{-T} \mecNref 	
	\quad \:\:\:\:\:\:\:\:\:\:\, \qquad \qquad
	\text{on }\Gamma^{LV}  \times (0, T], \label{eqn: Mec3}\\
	& \tenspiola(\displ_H, \Tens(\ActStateHF)) \mecNref = -\PRV(t) \, J_H \mecF_H^{-T} \mecNref
	\:\:\:\:\:\:\:\:\:\:\:\, \qquad \qquad
	\text{on }\Gamma^{RV} \times (0, T], \label{eqn: Mec4}\\
	\begin{split}& \tenspiola(\displ_H, \Tens(\ActStateHF)) \mecNref \\ & \quad = \displaystyle | J_H \mecF_H^{-T} \mecNref | \left[ \PLV(t)\BCmecVbaseLV + 
	\PRV(t)\BCmecVbaseRV \right] \qquad \:\:\:\, \text{on } \Gamma_H^{\text{base}} \times (0, T], \end{split} \label{eqn: Mec5}\\
	&\mathbf{d}_H= \mathbf{d}_{0,H}   \qquad \qquad \qquad \qquad \qquad \qquad \quad \:\:\:\:\, \text{in }\{\Omega_{H}^0 \cup \Omega_{C}^0\}\times \{t = 0\},
	\end{empheq}
\end{subequations}
where $\BCmecVbaseLV$ and $\BCmecVbaseRV$ are the following vectors \cite{REGAZZONI2022111083}:
$$ \mathbf{v}^{\text{base}}_{i} = \frac{\int_{\Gamma^{\text{endo},i}_H} J\mathbf{F}_H^{-T}\mathbf{N}d\Gamma_H^{i}}{\int_{\Gamma^{\text{endo},i}_H} |J\mathbf{F}_H^{-T}\mathbf{N}|d\Gamma_H^{i}} \qquad i = LV,RV.$$

The myocardial tissue is assumed to be an hyperelastic material \cite{guccione1991finite,ogden1997non}, while active mechanics is described by means of an active stress approach \cite{goktepe2010electromechanics,smith2004multiscale}. The Piola-Kirchhoff $\tenspiola = \tenspiola(\displ,\Tens)$ stress tensor is decomposed in a first term, representing the strain energy density function $\mathcal{W}:\mathrm{Lin}^+\to \mathbb{R}$, and a second one corresponding to the orthotropic active stress, namely:
$$ \tenspiola(\displ_H,\Tens) = \frac{\partial \mathcal{W}(\mecF_H)}{\partial \mecF_H} + \Tens(\hat{\xi},\ActStateHF)\left[ n_f \frac{\mecF_H\fZero \otimes \fZero}{\sqrt{\IIVf}} + n_n \frac{\mecF_H\nZero \otimes \nZero}{\sqrt{\IIVn}}   \right].$$
Here, $\mecF_H$ is the deformation tensor, while $\Tens(\hat{\xi},\ActStateHF)$ is the active tension provided by the activation model. $\IIVf = \mecF_H \fZero \cdot \mecF_H \fZero$ and $\IIVn = \mecF_H \nZero \cdot \mecF_H \nZero$ represent the tissue stretches along the fiber and sheet-normal directions, respectively, being $n_f$ and $n_n$ the prescribed portion of active stress tensor in fiber and sheet-normal directions.

The strain energy function $\mathcal{W}$ is described by the Guccione constitutive law \cite{guccione1991finite}:
\begin{equation}\label{Eq:guccione_law} \mathcal{W} = \frac{\kappa}{2}(J-1)\log(J) + \frac{\tilde{a}}{2}(e^Q-1),\end{equation}
where the first term accounts for the volumetric energy, including the bulk mudulus $\kappa$. The term $\tilde{a}$ is a stiffness scaling parameter $\tilde{a} = a[\mu+(1-\mu)4.56]$, see \cite{salvador2021electromechanical}, where $\mu$ is the parameter that account for possible scars and gray zones in the myocardium. 

In equation \eqref{Eq:guccione_law}, the exponent $Q$ is related to the Green-Lagrange strain tensor $\mathbf{E} = \frac{1}{2}(\mecC-\identity)$, being $\mecC = \mecF_H^T\mecF_H$ the right Cauchy-Green deformation tensor, by:
$$ Q = b_{\mathrm{ff}} E_{\mathrm{ff}}^2 + b_{\mathrm{ss}} E_{\mathrm{ss}}^2 + b_{\mathrm{fs}} (E_{\mathrm{fs}}^2+E_{\mathrm{sf}}^2)+b_{\mathrm{fn}} (E_{\mathrm{fn}}^2+E_{\mathrm{nf}}^2)+b_{\mathrm{sn}} (E_{\mathrm{sn}}^2+E_{\mathrm{nf}}^2), $$
where $b$ is the stiffness scaling parameter and $E_{\mathrm{ij}} = \mathbf{E}\mathbf{i}_0\cdot\mathbf{j}_0$, for i,j $\in \{$f,s,n$\}$ and $\mathbf{i}_0,\mathbf{j}_0 \in\{ \fZero,\nZero,\sZero\}$ are the entries of the Green-Lagrange strain tensor $\mathbf{E}$.

Boundary conditions \eqref{eqn: Mec2}-\eqref{eqn: Mec5} are prescribed to model the interaction of the endocardium with the blood, as well as the tension due to the continuity of the heart muscle on the base and the pericardium~\cite{Gerbi2018monolithic,pfaller2019importance,strocchi2020simulating}. The blood-endocardium interaction is modeled through normal stress boundary conditions \eqref{eqn: Mec3}-\eqref{eqn: Mec4} prescribed at the endocardial surface $\Gamma_H^{\text{endo,LV}}$ and $\Gamma_H^{\text{endo,RV}}$ and on the cap surfaces $\Gamma_C^{\text{endo,LV}}$ and $\Gamma_C^{\text{endo,RV}}$. The energy-consistent boundary condition accounting for the effect of the neglected part of the biventricular domain is instead imposed through \eqref{eqn: Mec5} on $\Gamma_H^{\text{base}}$ (see \cite{piersanti2021closedloop,regazzoni2020machine}) . Finally, the effect of the pericardium is accounted for by means of generalized Robin boundary conditions at the epicardial surfaces $\Gamma_H^{\text{epi}}$ \eqref{eqn: Mec2} trough the tensors $\mathbf{K}^{\mathrm{epi}} = K_{\|}^{\mathrm{epi}} (\mathbf{N}\otimes\mathbf{N}-\mathbf{I}) - K_{\perp}^{\mathrm{epi}}(\mathbf{N}\otimes\mathbf{N})$ and $\mathbf{C}^{\mathrm{epi}}= C_{\|}^{\mathrm{epi}}(\mathbf{N}\otimes\mathbf{N}-\mathbf{I}) - C_{\perp}^{\mathrm{epi}}(\mathbf{N}\otimes\mathbf{N})$, with $K_{\perp}^{\mathrm{epi}},~C_{\perp}^{\mathrm{epi}},~K_{\|}^{\mathrm{epi}},~C_{\|}^{\mathrm{epi}} \in \mathbb{R}^+$ the stiffness and viscosity parameters of the epicardium in normal and tangential directions, respectively.  

\subsection{Circulation}
The role of blood circulation in the cardiac contration is included by the 0D description of the complete cardiovaluscar system proposed in \cite{piersanti2021closedloop,REGAZZONI2022111083} that models the systemic and the pulmonary circulations as RLC circuits, the heart chambers with time-varying elanstance elements, and the heart valves through non-ideal diodes. The resulting ODE system:
\begin{subequations}
	\begin{empheq}[left={\empheqlbrace\,}]{align} \label{eqn: Circ}
	&\dfrac{d \Circ(t)}{d t} = \CircRhs(t, \Circ(t), \PLV(t), \PRV(t))   \qquad t \in (0,T],\\
	&\mathbf{c} = \mathbf{c}_0 \qquad \qquad \qquad \qquad \qquad \qquad \qquad \:\:\:\:\:\, t = 0,
	\end{empheq}
\end{subequations}
represents therefore the blood circulation. Pressures, volumes and fluxes of the different vascular compartments are included in the unknowns vector $\mathbf{c}$.

The coupling of the 0D circulatory model  with the 3D biventricular model EM model is achieved by replacing the time-varying elastance elements of the LV and RV with their corresponding 3D descriptions in the circulation model. A suitable volume-consistency coupling conditions, given by: 
$$V_i^{3D}(\mathbf{c}(t)) = \int_{\GammaEndoi}J(t)((\mathbf{h}\otimes\mathbf{h})(\mathbf{x}+\mathbf{d}_H(t)-\mathbf{b}_i))\cdot\mecF^{-T}(t)\mathbf{N}d\Gamma_0, \quad i=\textrm{LV, RV}$$
is introduced, being $\mathbf{h}$ an orthogonal vector to the LV and RV centerline, while $\mathbf{b}_i$ is a vector inside the LV and RV. Therefore, pressures of the LV and RV in the 3D-0D coupled model can be determined through the Lagrange multipliers associated to the constraints: 
\begin{subequations}
	\begin{empheq}[left={\empheqlbrace\,}]{align}
	\VLV(\Circ(t)) = \VLVthreedim(\displ_H(t))  & \qquad t \in (0,T], \label{eqn: Coupl1}\\
	\VRV(\Circ(t)) = \VRVthreedim(\displ_H(t))  &\qquad  t \in (0,T], \label{eqn: Coupl2}
	\end{empheq}
\end{subequations}
rather than via the 0D circulation model.

\subsection{12 lead ECG system}
The standard 12-lead ECG is a system of 12 leads obtained by combining the values of $u_T$ recorded from 9 electrodes on the surface of the human body, named $R$, $L$, $F$, and $V_i$, $i = 1,\dots,6$ (referred to Figure \ref{Fig.BSPM_healthy} for a representation of the electrodes distribution). Defining by $\mathbf{x}_R$, $\mathbf{x}_L$, $\mathbf{x}_F$, $\mathbf{x}_{V_i}$, the spatial location of the electrodes, the 6 limb leads are computed as:
\begin{gather*}
I = u_T(\mathbf{x}_L) - u_T(\mathbf{x}_R), \quad II = u_T(\mathbf{x}_F) - u_T(\mathbf{x}_R), \quad III = u_T(\mathbf{x}_F) - u_T(\mathbf{x}_L),\\ 
aVR = u_T(\mathbf{x}_R) - \frac{1}{2} (u_T(\mathbf{x}_L) + u_T(\mathbf{x}_F)), \quad aVL = u_T(\mathbf{x}_L) - \frac{1}{2}(u_T(\mathbf{x}_R) + u_T(\mathbf{x}_F)), \\ aVF = u_T(\mathbf{x}_F) - \frac{1}{2}(u_T(\mathbf{x}_L) + u_T(\mathbf{x}_R)),
\end{gather*}
whereas the 6 precordial (or chest) leads are defined as:
\[ V_i = u_T(\mathbf{x}_{V_{i}}) -WCT, \quad i = 1,\dots,6,
\]
with $WCT$ denoting the Wilson central terminal signal, given by:
$$WCT =  \frac{1}{3} \left[ u_T(\mathbf{x}_L) + u_T(\mathbf{x}_R) + u_T(\mathbf{x}_F)\right ].$$

%
%
%

	\end{appendix}
\end{document}